\documentclass[12pt,a4paper,oneside]{amsart}
\usepackage{mathrsfs}
\usepackage{amssymb}
\usepackage{amsxtra}
\usepackage{hyperref, cleveref}
\usepackage{graphics,color,graphicx}
\usepackage[T1]{fontenc}
\usepackage[utf8]{inputenc}
\usepackage{tikz}
\usepackage{caption}
\usepackage{xcolor}
\usepackage{marginnote} 

\newcommand{\fC}{\mathfrak{C}}
\newcommand{\fD}{\mathfrak{D}}

\newcommand{\fF}{\mathfrak{F}}
\newcommand{\fG}{\mathfrak{G}}

\newcommand{\fK}{\mathfrak{K}}
\newcommand{\fL}{\mathfrak{L}}
\newcommand{\fM}{\mathfrak{M}}
\newcommand{\fN}{\mathfrak{N}}

\newcommand{\fT}{\mathfrak{T}}

\newcommand{\pA}{\mathcal{A}}
\newcommand{\pB}{\mathcal{B}}

\newcommand{\pD}{\mathcal{D}}

\newcommand{\pL}{\mathcal{L}}

\newcommand{\pO}{\mathcal{O}}

\newcommand{\pS}{\mathcal{S}}
\newcommand{\pT}{\mathcal{T}}

\newcommand{\eH}{\mathscr{H}}

\newcommand{\eK}{\mathscr{K}}

\newcommand{\eM}{\mathscr{M}}
\newcommand{\eN}{\mathscr{N}}

\newcommand{\eU}{\mathscr{U}}
\newcommand{\eV}{\mathscr{V}}
\newcommand{\eX}{\mathscr{X}}
\newcommand{\eW}{\mathscr{W}}

\newcommand{\bC}{\mathbb{C}}

\newcommand{\bR}{\mathbb{R}}

\newcommand{\bN}{\mathbb{N}}
\newcommand{\bQ}{\mathbb{Q}}

\newcommand{\bZ}{\mathbb{Z}}

\newcommand{\varBeta}{B}
\newcommand{\sot}{{\it sot}}

\DeclareMathOperator{\Ker}{{\rm ker}}

\DeclareMathOperator{\Alg}{Alg}

\DeclareMathOperator{\Aut}{Aut}
\DeclareMathOperator{\Col}{Col}

\DeclareMathOperator{\Lat}{Lat}

\DeclareMathOperator{\Grp}{{\rm Grp}}

\DeclareMathOperator{\supp}{{\rm supp}}

\newtheorem{theorem}{Theorem}
\newtheorem{proposition}[theorem]{Proposition}
\newtheorem{lemma}[theorem]{Lemma}
\newtheorem{corollary}[theorem]{Corollary}
\theoremstyle{definition}
\newtheorem{definition}[theorem]{Definition}

\newtheorem{example}[theorem]{Example}

\numberwithin{equation}{section}
\begin{document}
%
\title[Similarities of subspace lattices in Banach spaces]{Similarities of subspace lattices in Banach spaces}
\author[J. Bra\v{c}i\v{c}]{Janko Bra\v{c}i\v{c}}
\address{Faculty of Natural Sciences and Engineering, University of Ljubljana, 
A\v{s}ker\v{c}eva c. 12, SI-1000 Ljubljana, Slovenia}
\email{janko.bracic@ntf.uni-lj.si}
\author{Marko Kandi\'c}
\address{Faculty of Mathematics and Physics, University of Ljubljana,
  Jadranska 19, 1000 Ljubljana, Slovenia \ \ \ and \ \ \
  Institute of Mathematics, Physics, and Mechanics, Jadranska 19, 1000 Ljubljana, Slovenia}
\email{marko.kandic@fmf.uni-lj.si}
%
\begin{abstract}
A collineation of a subspace lattice $\fL$ in a complex Banach space $\eX$ is an invertible operator $S$ on $\eX$ 
with the property that the image $S\eM$ of a subspace $\eM$ belongs to $\fL$ if and and only if $\eM$ belongs to it.
Hence, $S$ is a collineation of $\fL$ if and only if it implements an order automorphism of $\fL$. We study the group 
$\Col(\fL)$ of all collineations of $\fL$ and its subgroup $\Grp(\Alg(\fL))$ of all invertible operators that fix 
every subspace in $\fL$. We show that $\Grp(\Alg(\fL))$ is a normal subgroup of $\Col(\fL)$; moreover, if $\fL$
is a reflexive subspace lattice, then $\Col(\fL)$ is the normalizer of $\Grp(\Alg(\fL))$ in the group of all invertible
operators on $\eX$. One of the main questions that we consider is whether $\Grp(\Alg(\fL))$ is a complemented
subgroup in $\Col(\fL)$. For certain subspace lattices $\fL$, such as some realizations of the diamond or the 
double triangle, some nests in the space of continuous functions on $[0,1]$, and the classical Volterra nest in 
$L^1[0,1]$, we characterize the complement of $\Grp(\Alg(\fL))$ in $\Col(\fL)$. On the other hand, for the Volterra
nests in $L^p[0,1]$, where $1<p<\infty$, a further study is needed, and we prove only some partial results.
\end{abstract}
\maketitle
%
\section{Introduction} \label{Sec01}
\setcounter{theorem}{0}

An invertible operator $S$ on a complex Banach space $\eX$ is a collineation of a non-empty family $\fF$ of
subspaces of $\eX$ if $\eM\in \fF\; \iff\; S\eM\in \fF$. When $\fF$ is a lattice of subspaces, collineations are 
referred to as similarities of the subspace lattice by some authors. However, to distinguish the notion from other types
of similarities, we keep the terminology that we introduced in \cite{BK}. The main object of our study in this paper
is the group $\Col(\fF)$ of all collineations of $\fF$. This study has been initiated in \cite{BK} in the context of
finite-dimensional complex vector spaces. Here we continue with our research and extend it to arbitrary complex
Banach spaces.

Let $\eX$ be a complex Banach space, $\pB(\eX)$ the algebra of all bounded linear operators on $\eX$,
and $\fC(\eX)$ the lattice of all closed subspaces of $\eX$. For a unital algebra $\pA\subseteq\pB(\eX)$, let
$\Grp(\pA)$ be the group of all operators that are invertible in $\pA$. For a subspace lattice $\fL\subseteq \fC(\eX)$,
group $\Grp\bigl(\Alg(\fL)\bigr)$ is a normal subgroup of the group of collineations $\Col(\fL)$. Moreover,
if $\fL$ is reflexive in the sense that $\fL=\Lat\Alg(\fL)$, then $\Col(\fL)$ is the normalizer of
$\Grp\bigl(\Alg(\fL)\bigr)$ in the group of all invertible operators on $\eX$ (\Cref{theo01}). Collineations
of a reflexive lattice $\fL$ can be characterized as those invertible operators on $\eX$ which implement
spatial automorphisms of $\Alg(\fL)$ or, equivalently, those operators which preserve the family of
cyclic subspaces of $\Alg(\fL)$ (\Cref{theo02}). Since $\Grp\bigl(\Alg(\fL)\bigr)$ is a normal subgroup
of $\Col(\fL)$ whenever $\fL$ is a subspace lattice, it is a natural question whether $\Grp\bigl(\Alg(\fL)\bigr)$
is complemented in $\Col(\fL)$, that is, does there exist a subgroup $\pO(\fL)\subseteq \Col(\fL)$ such that
$\Col(\fL)$ is a semidirect product of $\Grp\bigl(\Alg(\fL)\bigr)$ and $\pO(\fL)$. The answer is obvious if the 
group of order automorphisms of $\fL$ is trivial (\Cref{prop05}), as in this case $\Col(\fL)$ must be equal to 
$\Grp\bigl(\Alg(\fL)\bigr)$. For medial lattices and nests, the situation is more interesting. For subspace lattices 
$\fD$ and $\fT$, which are realizations of the diamond and the double triangle, respectively, by isomorphic
algebraically complemented subspaces, one has $\Col(\fD) = \Grp\bigl(\Alg(\fD)\bigr) \rtimes \pO(\fD)$ and 
$\Col(\fT) = \Grp\bigl(\Alg(\fT)\bigr) \rtimes \pO(\fT)$, where $\pO(\fD) \subseteq \Col(\fD)$ is a subgroup
isomorphic to the symmetric group $S_2$, and $\pO(\fT) \subseteq \Col(\fT)$ is a subgroup isomorphic to the
symmetric group $S_3$ (\Cref{theo04}).
If $\fN=\{ \eN_k;\; k\in \bZ\}\cup\bigl\{ \{0\},\ell^p(\bZ)\bigr\}$ is the classical nest of subspaces of $\ell^p(\bZ)$
$(1\leq p<\infty)$ where $\eN_k=\bigvee\{ e_j;\; j\leq k\}$, then 
$\Col(\fN)=\Grp\bigl(\Alg(\fN)\bigr)\rtimes\pO(\fN)$, where $\pO(\fN)$ is the group of isometries generated by 
the unilateral shift (\Cref{prop12}). A similar result holds for the Volterra nest of subspaces of $L^1[0,1]$ but not 
for the Volterra nest of subspaces of $L^2[0,1]$ (\Cref{theo08}).

\section{Preliminaries} \label{Sec02}
\setcounter{theorem}{0}

Since the main objects of our study will be subspace lattices and their similarities, we begin with a few general facts
from lattice theory (for more details see \cite{Bir,Gra}).
By a lattice we mean an algebraic structure $(L,\vee,\wedge)$, consisting of a non-empty set $L$
and two binary operations (meet and join) that are commutative, associative and satisfy the absorption laws
$a\vee(a\wedge b)=a$ and $a\wedge(a\vee b)=a$ for all $a, b\in L$.
One can define a partial order $\leq$ on $L$ by setting $a\leq b$ if $a=a\wedge b$ or, equivalently, $b=a\vee b$.
A lattice $L$ is bounded if it has a greatest element $1$ and a least element $0$. Hence, $0\leq a\leq 1$,
for all $a\in L$. In this paper, we work only with bounded lattices.
An element $a\in L$ is non-trivial if $a\ne 0$ and $a\ne 1$.
Elements $a, b\in L$ are comparable if either $a\leq b$ or $b\leq a$, otherwise they are incomparable.
An element $a\in L$ is said to be comparable if it is comparable with every element in $L$.

Examples of lattices in which we are interested are the family $\fC(\eX)$ of all closed subspaces of a complex
Banach space $\eX$ and its sublattices. It is worth mentioning that the meet and the join operations in
$\fC(\eX)$ are defined by $\eM\wedge\eN=\eM\cap\eN$ and $\eM\vee\eN=\overline{\eM+\eN}$
for arbitrary subspaces $\eM, \eN\in \fC(\eX)$. Hence, the partial order is the set inclusion.
Of course, if $\eX$ is finite-dimensional, then $\fC(\eX)$ is the lattice of all linear submanifolds of $\eX$.
A family $\fL\subseteq\fC(\eX)$ that contains the trivial subspaces $\{0\}$ and $\eX$ will be called
a {\em lattice of subspaces} if it is a sublattice of $\fC(\eX)$, that is, $\eM\wedge\eN$ and $\eM\vee\eN$ are in
$\fL$ whenever $\eM$ and $\eN$ are in $\fL$.

An isomorphism of lattices $L_1$ and $L_2$ is a bijective mapping $\phi\colon L_1\to L_2$ such that
$\phi(a\vee b)=\phi(a)\vee\phi(b)$ and $\phi(a\wedge b)=\phi(a)\wedge\phi(b)$ for all $a, b\in L_1$.
It is not hard to see that a bijective mapping $\phi\colon L_1\to L_2$ is an isomorphism if and only if it
is an order isomorphism, that is, $a\leq b\iff \phi(a)\leq \phi(b)$ for all $a, b\in L_1$.
An automorphism of a lattice $L$ means an isomorphism $\phi\colon L\to L$.
We will use $\Aut(L)$ to denote the group of all automorphisms of a lattice $L$.
A bijective mapping $\psi\colon L_1\to L_2$ such that $a\leq b\iff \psi(a)\geq \psi(b)$ (equivalently,
$\psi(a\vee b)=\psi(a)\wedge\psi(b)$ and $\psi(a\wedge b)=\psi(a)\vee\psi(b)$) for all $a, b\in L_1$,
is an anti-isomorphism.

In the case of lattices of subspaces, isomorphisms $\Phi\colon\fC(\eX)\to\fC(\eX)$ are sometimes called projectivities
of $\eX$ (see \cite[p.\ 40]{Bae}). If $\eX$ is finite-dimensional and $\dim(\eX)\geq 3$, then
the first fundamental theorem of projective geometry says that every projectivity of $\eX$ is implemented by a
semi-linear transformation. This means that for every projectivity $\Phi$ of $\eX$ there exists an invertible semi-linear
mapping $S\colon\eX\to\eX$ such that $\Phi(\eM)=S\eM$ for all $\eM\in \fC(\eX)$.
Recall that $S$ is semi-linear if it is additive and $S(\lambda x)=\tau(\lambda) Sx$ for all $\lambda\in\bC$ and
all $x\in \eX$, where $\tau$ is an automorphism of $\bC$. A projectivity of $\eX$ that is implemented by an
invertible semi-linear mapping $S$ will be denoted by $\Phi_S$. Fillmore and Longstaff \cite{FL} proved that every
projectivity of an infinite-dimensional complex Banach space $\eX$ is implemented by a bounded linear or
a conjugate linear invertible operator. Every invertible linear
or conjugate linear operator $S$ clearly induces a projectivity $\Phi_S\colon \eM\mapsto S\eM$ of $\eX$.

A lattice $L$ with ordering $\leq$ is complete if every subset $A\subseteq L$ has both a greatest lower bound
$\bigwedge A$ and a least upper bound $\bigvee A$ in $L$. Notice that $\bigwedge\emptyset=1$ and
$\bigvee\emptyset=0$. If $\fF$ is a family of subspaces in $\fC(\eX)$, then $\bigwedge\fF$ is the intersection of
all subspaces in $\fF$ and $\bigvee \fF$ is the closed linear span of all subspaces in $\fF$. In particular,
$\bigwedge\emptyset=\eX$ and $\bigvee\emptyset=\{0\}$. A complete lattice of subspaces will be called a
{\em subspace lattice}.

We will frequently use a known fact that isomorphisms between complete lattices are complete.
More precisely, if $L_1$ and $L_2$ are complete lattices and $\phi\colon L_1\to L_2$ is an order isomorphism, then
$\phi\left(\bigvee\{ a_j; j\in J\}\right)=\bigvee\{ \phi(a_j); j\in J\}$ and
$\phi\left(\bigwedge\{ a_j; j\in J\}\right)=\bigwedge\{ \phi(a_j); j\in J\}$ for an arbitrary family
$\{ a_j;\; j\in J\}\subseteq L_1$.

A subspace $\eM\in \fC(\eX)$ is invariant for $T\in \pB(\eX)$ if $T\eM\subseteq \eM$.
More generally, for a non-empty set $\pT\subseteq \pB(\eX)$, a subspace $\eM$ is invariant for $\pT$ if it is
invariant for every operator in $\pT$. Let $\Lat(\pT)\subseteq \fC(\eX)$ denote the subspace
lattice of all subspaces that are invariant for $\pT$. We denote the commutant of $\pT$ by $\pT'$. 
It is not hard to see that $\pT'$ is a unital \sot-closed subalgebra of $\pB(\eX)$ (unital means that it contains $I$, 
the identity operator on $\eX$, and \sot\ is an abbreviation for the strong operator topology). 
The bicommutant of $\pT$ is $\pT''=(\pT')'$. It is easily seen that $\pT\subseteq \pT''$ and $\pT'''=\pT'$.

For a non-empty family $\fF\subseteq \fC(\eX)$, let
$\Alg(\fF)=\{ T\in \pB(\eX);\; T\eM\subseteq \eM,\;\text{for all}\; \eM\in \fF\}$. It is a unital \sot-closed
subalgebra of $\pB(\eX)$. The reflexive cover of a non-empty set $\pT\subseteq \pB(\eX)$ is $\Alg\Lat(\pT)$.
A \sot-closed subalgebra $\pA\subseteq \pB(\eX)$ is reflexive if  $\Alg\Lat(\pA)=\pA$. Notice that
this happens if and only if there exists a non-empty family $\fF\subseteq \fC(\eX)$ such that $\pA=\Alg(\fF)$.
A subspace lattice $\fL\subseteq \fC(\eX)$ is reflexive if $\fL=\Lat\Alg(\fL)$.

Let $\pA\subseteq \pB(\eX)$ be a unital algebra. If $S \in \pA$ is an invertible operator, its inverse $S^{-1}$ does not 
necessarily belong to $\pA$. By $\Grp(\pA)$ we will denote the group of all invertible
operators $S\in \pA$ such that $S^{-1}\in\pA$. In particular, $\Grp\bigl(\pB(\eX)\bigr)$ is the group of all invertible
operators on $\eX$. If $\fF\subseteq\fC(\eX)$ is a non-empty family, then an invertible operator $T$ belongs to
$\Grp\bigl(\Alg(\fF)\bigr)$ if and only if $T\eM=\eM$ for all $\eM\in \fF$. Indeed, if $T\in \Grp\bigl(\Alg(\fF)\bigr)$,
then $T, T^{-1}\in \Alg(\fF)$ and therefore $T\eM\subseteq \eM$ and $T^{-1}\eM\subseteq \eM$ for all
$\eM\in \fF$. It follows that $T\eM=\eM$ for all $\eM\in \fF$. For the opposite implication, observe that
$T\eM=\eM$ gives $T^{-1}\eM=\eM$ for every $\eM\in \fF$. Hence, $T, T^{-1}\in\Alg(\fF)$.

\section{Collineations} \label{Sec03}
\setcounter{theorem}{0}

\subsection{Definition and basic observations}  \label{Sec031}

Subspace lattices $\fK$ and $\fL$ of $\fC(\eX)$ are similar if there exists an invertible operator $S\in\pB(\eX)$
such that the mapping $\Phi_S\colon \eM\to S\eM$ is an order isomorphism from $\fK$ onto $\fL$
(see \cite{AAW, ALWW, Dav1, Dav, Lar, Lar1}). In this case, $S$ is called a similarity between
$\fK$ and $\fL$. On the other hand, in the context of projective geometry, projectivities of $\eX$ that are implemented
by linear operators are usually called collineations (see \cite[p. 62]{Bae}).
With a slight abuse of terminology, we introduce the following definition (cf. \cite{BK}).

\begin{definition} \label{def01}
An invertible operator $S\in \pB(\eX)$ is a {\em collineation} of a non-empty family of subspaces
$\fF\subseteq \fC(\eX)$ if $\eM\in \fC(\eX)$ belongs to $\fF$ if and only if $S\eM$ belongs to $\fF$.
\end{definition}

We will denote the set of collineations of $\fF$ by $\Col(\fF)$. For any non-empty set $\pT$ of $\pB(\eX)$ we
define $\Col(\pT):=\Col\bigl(\Lat(\pT)\bigr)$ and we will call operators from $\Col(\pT)$ as collineations of $\pT$.
Since $\Lat(\pT) = \Lat\Alg\Lat(\pT)$, there is no loss of generality in focusing on reflexive subalgebras of 
$\pB(\eX)$.

\begin{lemma} \label{lem03}
Let $\fF\subseteq \fC(\eX)$ be a non-empty family. The following assertions are equivalent for an invertible
operator $S$.
\begin{itemize}
\item[(a)] $S$ is a collineation of $\fF$;
\item[(b)] if $\eM\in\fF$, then $S\eM$ and $S^{-1}\eM$ are in $\fF$;
\item[(c)] $S^{-1}$ is a collineation of $\fF$.
\end{itemize}

\noindent
Hence, $\Col(\fF)$ is a group containing the group $\Grp\bigl(\Alg(\fF)\bigr)$.
If $S_1, S_2\in\Col(\fF)$, then $S_1\eM=S_2\eM$ for all $\eM\in\fF$ if and only if
$S_{2}^{-1}S_1\in\Grp\bigl(\Alg(\fF)\bigr)$.
\end{lemma}

\begin{proof}
We will check only the equivalence of (a) and (b) since it implies the
equivalence of (c) and (b). Assume (a) and let $S$ be a collineation of $\fF$. If $\eM \in \fL$, then, trivially, 
$S\eM \in \fF$, and since $S(S^{-1}\eM) = \eM \in \fF$, it follows that $S^{-1}\eM \in \fF$. Suppose that (b) holds.
We only have to see that $S\eM\in \fF$ implies $\eM\in \fF$ which follows from $\eM=S^{-1}(S\eM)\in \fF$.

It is obvious now that $\Col(\fF)$ is a group. If $S\in \Grp(\Alg(\fF))$, then $S\eM=\eM\in \fF$ and 
$S^{-1}\eM=\eM\in \fF$ for every $\eM\in \fF$, proving $\Grp\bigl(\Alg(\fF)\bigr)\subseteq \Col(\fF)$.
The last statement is obvious.
\end{proof}

\Cref{lem03} is very simple; we have explicitly stated it to stress that $S$ is a collineation of $\fF$ if both $S$
and $S^{-1}$ preserve $\fF$. The following example shows that the condition
\begin{equation} \label{eq02}
\eM\in\fF\quad\Longrightarrow\quad S\eM\in \fF
\end{equation}
is usually not enough for an invertible operator $S$ to be a collineation of $\fF$. 

\begin{example} \label{ex04}
Consider the Hilbert space $\ell^2(\mathbb Z)$ of all square-summable two-sided sequences with the standard
orthonormal basis $\{e_n;\; n\in \mathbb Z\}$. Let $W$ be the bilateral shift on $\ell^2(\mathbb Z)$ given by 
$W e_n=e_{n+1}$ for all $n\in \bZ$. For every $n\in \bZ$, let $\eM_n=\bigvee\{ e_j;\; j\geq n\}$ and 
$\fF=\{ \eM_n;\; n\geq 0\}$. Since $W\eM_n=\eM_{n+1}\subsetneq \eM_n$ for all $n\in \bZ$, the operator 
$W$ satisfies \eqref{eq02}. On the other hand,  $W^{-1}\eM_0=\eM_{-1} \notin \fF$ which means that 
$W$ is not a collineation of $\fF$.\qed
\end{example}

Nevertheless, for some non-empty families $\fF$, \eqref{eq02} implies $S^{-1}\eM\in \fF$ for every $\eM\in \fF$. 
For instance, this is the case when $\eX$ is a finite-dimensional Banach space and $\fF$ is  a reflexive subspace 
lattice (see \cite[Corollary 2.3]{BK}). \Cref{ex04} shows that in infinite-dimensional Banach spaces condition
\eqref{eq02} does not guarantee that $S$ is a collineation of $\fF$ even when $\fF$ is a reflexive subspace lattice.
Indeed, an easy verification shows that $\fF$ is a totally ordered subspace lattice which is reflexive by 
\cite[Theorem 3.4]{Rin1}. 

The following proposition provides two examples in infinite-dimensional Banach spaces for which 
\eqref{eq02}characterizes $\Col(\fF)$. Recall that a subspace lattice is a nest if it is totally ordered. 
A nest is maximal if it is not contained in any larger nest.

\begin{proposition} \label{prop15}
If $\fF\subseteq \fC(\eX)$ is a non-empty finite family or a maximal nest, then an invertible operator $S$ is
a collineation of $\fF$ if and only if \eqref{eq02} holds.
\end{proposition}

\begin{proof}
If $S\in\Col(\fF)$, then \eqref{eq02} holds, by the definition of a collineation. We have to prove the opposite
implication. Suppose first that $\fF=\{\eM_1,\ldots,\eM_k\}$, where $\eM_j$ ($j=1,\ldots,k$) are distinct subspaces 
of $\eX$. Let $S\in \pB(\eX)$ be an invertible operator such that $S\eM_j\in \fF$
for every $j=1,\ldots, k$. Then $S\eM_j$ ($j=1,\ldots,k$) are distinct subspaces in $\fF$, and therefore,
there exists a permutation $\pi$ of $\{1,\ldots, k\}$ such that $S\eM_j=\eM_{\pi(j)}$. It follows that
$S^{-1}\eM_j=\eM_{\pi^{-1}(j)}$ $(j=1,\ldots, k)$. Hence, by \Cref{lem03}, we have $S\in \Col(\fF)$.

Assume that $\fF$ is a maximal nest. If $S$ is an invertible operator satisfying \eqref{eq02}, then,
$\eM\subseteq S\eN$ or $S\eN\subseteq \eM$ for arbitrary $\eM, \eN\in \fF$. Hence, for a fixed $\eM$,
we have $S^{-1}\eM\subseteq \eN$ or $\eN\subseteq S^{-1}\eM$ for every $\eN\in \fF$. It follows,
by the maximality of $\fF$, that $S^{-1}\eM\in \fF$. Now we conclude, by \Cref{lem03}, that $S\in \Col(\fF)$.
\end{proof}

Given a non-empty family of subspaces $\fF \subseteq \fC(\eX)$, we denote by $\fL_{\fF} \subseteq \fC(\eX)$ the 
smallest lattice of subspaces that contains $\fF$. It is straightforward to verify that $\fL_{\fF}$ coincides with 
the intersection of all lattices of subspaces in $\fC(\eX)$ that contain $\fF$. Similarly, let 
$\widehat{\fL} \subseteq \fC(\eX)$ denote the smallest subspace lattice containing a given lattice of subspaces $\fL$. 
Then $\widehat{\fL}$ is precisely the intersection of all subspace lattices that contain $\fL$.

\begin{theorem} \label{theo07}
Let $\fF\subseteq \fC(\eX)$ be a non-empty family and let $\fL\subseteq \fC(\eX)$ be a lattice of subspaces.
If $\fF\subseteq \fG\subseteq \fL_\fF$, then $\Col(\fG)$ is a subgroup of $\Col(\fL_{\fF})$. Similarly, if
$\fK\subseteq\fC(\eX)$ is a lattice of subspaces such that $\fL\subseteq\fK\subseteq \widehat{\fL}$, then
$\Col(\fK)$ is a subgroup of $\Col(\widehat{\fL})$.
\end{theorem}

\begin{proof}
Pick $S\in\Col(\fG)$ and let us define $\fL':=\{\eM\in \fL_{\fF};\; S\eM\in \fL_{\fF}\}$. By assumption, $\fL'$ 
contains $\fG$. We claim that $\fL'$ is a lattice of subspaces contained in $\fL_{\fF}$. To this end, choose arbitrary
subspaces $\eM$ and $\eN$ in $\fL'$. Since $S$ is invertible, we have $S(\eM\vee\eN)=S\eM\vee S\eN$ and
$S(\eM\wedge\eN)=S\eM\wedge S\eN$. Applying the fact that $S\eM$ and $S\eN$ belong to the lattice of subspaces
$\fL_{\fF}$ we conclude that $\fL'$ is a lattice of subspaces. Since $\fL'$ contains $\fF$ which generates $\fL_{\fF}$,
we have $\fL'=\fL_{\fF}$. The same argument shows that $S^{-1}$ maps $\fL_{\fF}$ into $\fL_{\fF}$, so that by
\Cref{lem03} we conclude that  $\Col(\fG)$ is contained in $\Col(\fL_{\fF})$. That $\Col(\fG)$ is a subgroup of 
$\Col(\fL_{\fF})$ should be clear.  

To prove that $\Col(\fK)$ is a subgroup of $\Col(\widehat{\fL})$ we proceed similarly as in the first part of the proof. 
We first pick $S\in\Col(\fK)$ and then define $\fL':=\{\eM\in \widehat{\fL};\; S\eM\in \widehat{\fL}\}$.
Then $\fL'$ is a lattice of subspaces, which by assumption contains $\fK$. Since $S$ is invertible, for any family 
$\{\eM_j;\; j\in J\}\subseteq \fL'$, we have $S\bigl( \bigvee_{j\in J}\eM_j\bigr)=\bigvee_{j\in J}S\eM_j$
and $S\bigl( \bigwedge_{j\in J}\eM_j\bigr)= \bigwedge_{j\in J}S\eM_j$. Hence, the completeness of $\widehat{\fL}$ 
gives the completeness of $\fL'$ and therefore $\fL'$ is a subspace lattice that contains $\fL$. We conclude that
$\fL'=\widehat{\fL}$. We omit the rest of the proof as it is the same as the first part of the proof.
\end{proof}

In the end of \Cref{Sec052} we shall see that $\Col(\fL_{\fF})$ can be a proper subgroup of
$\Col(\widehat{\fL_{\fF}})$.

\subsection{Collineations of reflexive algebras}  \label{Sec032}

An automorphism of a subalgebra $\pA\subseteq \pB(\eX)$ is a bijective linear mapping $\Theta\colon\pA\to\pA$
such that $\Theta(AB)=\Theta(A)\Theta(B)$ for all $A,B\in\pA$. If there exists an invertible
operator $S\in\pB(\eX)$ such that $\Theta(A)=SAS^{-1}$ for all $A\in\pA$, then $\Theta$ is a
spatial automorphism. The spatial automorphism implemented by $S$ is denoted by $\Theta_S$. It should be clear
that $\Theta_S\colon \pA\to\pA$ is a spatial automorphism if and only if $S\pA S^{-1}=\pA$. Notice that the inverse
of a spatial automorphism $\Theta_S$ is again a spatial automorphism, as $\Theta_S^{-1}=\Theta_{S^{-1}}$.
It is easily seen that invertible operators $S_1, S_2\in \pB(\eX)$ induce the same spatial automorphism of a subalgebra
$\pA\subseteq \pB(\eX)$ if and only if $S_{2}^{-1}S_1\in \pA'$. A spatial automorphism of $\pA$ is inner if it
is implemented by an operator $S\in \Grp(\pA)$.

For a unital algebra $\pA\subseteq \pB(\eX)$ and a vector $x\in \eX$, let $\pA_x$ be the closure of the orbit
$\{ Ax;\; x\in \pA\}$. The subspace $\pA_x$ is called the cyclic subspace for $\pA$ generated by $x$.
It is clear that $\pA_x\in \Lat(\pA)$ for every $x\in \eX$, and it is easily seen that $\Lat(\pA)$ is generated by
the family of all cyclic subspaces. Hence, an operator $T\in \pB(\eX)$ belongs to $\Alg\Lat(\pA)$ if and only if
$T\pA_x\subseteq \pA_x$ for all $x\in \eX$.

The following lemma is known; it essentially says that every spatial automorphism of a unital algebra of operators 
can be extended to a spatial automorphism of the reflexive cover of the algebra.

\begin{lemma}\label{lem06}
Let $\pA\subseteq \pB(\eX)$ be a unital algebra. If $S \in \pB(\eX)$ is invertible and satisfies $S\pA S^{-1} = \pA$,
then $S\pA_x = \pA_{Sx}$ for every $x \in \eX$ and $S\Alg\Lat(\pA)S^{-1} = \Alg\Lat(\pA)$.
\end{lemma}

\begin{proof}
Let $x\in \eX$ be arbitrary. Since $S^{-1}\pA S=\pA$ we have $\pA_x=\overline{\{ S^{-1}ASx;\, A\in \pA\}}$
and therefore 
$$S\pA_x=S\overline{\{ S^{-1}ASx;\, A\in \pA\}}=\overline{S\{S^{-1}ASx;\, A\in \pA\}}=
\overline{\{ASx;\, A\in \pA\}}=\pA_{Sx}.$$ 
Notice that we have $S^{-1}\pA_x=\pA_{S^{-1}x}$ for every $x\in\eX$,
as well. Hence, if $T\in \Alg\Lat(\pA)$, then $STS^{-1}\pA_x=ST\pA_{S^{-1}x}\subseteq S\pA_{S^{-1}x}=\pA_x$
for every $x\in \eX$. Since cyclic subspaces generate $\Lat(\pA)$, it follows that
$S\Alg\Lat(\pA) S^{-1}\subseteq\Alg\Lat(\pA)$. Similarly, $S^{-1}\Alg\Lat(\pA) S\subseteq\Alg\Lat(\pA)$.
\end{proof}

Recall that the normalizer of a subgroup $H$ of a group $G$ is the subgroup
$N_{_G}(H)=\{ g\in G;\; gHg^{-1}=H\}$ of $G$. If $N_{_G}(H)=G$, then $H$ is said to be a normal subgroup
of $G$. By definition, the normalizer $N_{_G}(H)$ of $H$ in $G$ is the largest subgroup of $G$ in which $H$ is normal.

\begin{theorem} \label{theo01}
For a subspace lattice $\fL\subseteq \fC(\eX)$, the following assertions hold.
\begin{itemize}
\item[(i)] The group $\Grp\bigl(\Alg(\fL)\bigr)$ is a normal subgroup of $\Col(\fL)$.
\item[(ii)] The group $\Col(\fL)$ is a subgroup of $\Col\bigl(\Lat\Alg(\fL)\bigr)$.
\item[(iii)] If $\fL$ is reflexive, then $\Col(\fL)$ is the normalizer of $\Grp\bigl(\Alg(\fL)\bigr)$
in $\Grp\bigl(\pB(\eX)\bigr)$.
\item[(iv)] If $\Col(\fL)$ is the normalizer of $\Grp\bigl(\Alg(\fL)\bigr)$ in $\Grp\bigl(\pB(\eX)\bigr)$,
then $\Col\bigl(\Lat\Alg(\fL)\bigr)$ $=\Col(\fL)$.
\end{itemize}
\end{theorem}

\begin{proof}
(i) By \Cref{lem03}, $\Grp\bigl(\Alg(\fL)\bigr)$ is a subgroup of $\Col(\fL)$.
Let $C\in \Col(\fL)$ and $S\in \Grp\bigl(\Alg(\fL)\bigr)$ be arbitrary. If $\eM\in \fL$,
then $C^{-1}\eM\in \fL$ and therefore $SC^{-1}\eM=C^{-1}\eM$. It follows that
$CSC^{-1}\eM=CC^{-1}\eM=\eM$. Similarly,  $CS^{-1}C^{-1}\eM=\eM$. Hence,
$CSC^{-1}\in \Grp\bigl(\Alg(\fL)\bigr)$.

(ii) Let $S\in \Col(\fL)$ and $\eM\in \fL$ be arbitrary. Since $S\eM\in \fL$, it follows that $AS\eM\subseteq S\eM$ 
for every $A \in \Alg(\fL)$, which implies that $S^{-1}\Alg(\fL)S\subseteq \Alg(\fL)$. Since $S^{-1}\in \Col(\fL)$,
by replacing $S$ with $S^{-1}$ in the latter set inclusion we get $S^{-1}\Alg(\fL)S=\Alg(\fL)$.  

Let now $\eN\in \Lat\Alg(\fL)$ be arbitrary. Since $S^{-1}\Alg(\fL)S= \Alg(\fL)=\Alg\Lat\Alg(\fL)$, 
for every $A\in\Alg(\fL)$ we have $S^{-1}AS\eN\subseteq \eN$, and therefore, $AS\eN\subseteq S\eN$.
Hence, $S\eN\in\Lat\Alg(\fL)$. Similarly, $S^{-1}\eN\in\Lat\Alg(\fL)$. We conclude, by \Cref{lem03},
that $S\in\Col\bigl(\Lat\Alg(\fL)\bigr)$.

(iii) Suppose that $\fL$ is reflexive.  By (i), it suffices to prove that the normalizer of $\Grp\bigl(\Alg(\fL)\bigr)$ in
$\Grp\bigl(\pB(\eX)\bigr)$ is contained in $\Col(\fL)$. To this end, let $T\in \pB(\eX)$ be an invertible operator in the
normalizer of $\Grp\bigl(\Alg(\fL)\bigr)$. Then
$T^{-1}\Grp\bigl(\Alg(\fL)\bigr)T=\Grp\bigl(\Alg(\fL)\bigr)$. Hence, for $S\in \Grp\bigl(\Alg(\fL)\bigr)$ and
$\eM\in \fL$ we have $T^{-1}ST\eM=\eM$, and therefore, $ST\eM=T\eM$.
Hence, $T\eM\in \Lat\bigl(\Grp\bigl(\Alg(\fL)\bigr)\bigr)=\Lat\Alg(\fL)$, and since $\fL$ is reflexive, we conclude 
that $T\eM\in \Lat\Alg(\fL)=\fL$.
Similarly, $T^{-1}\eM\in \Lat\bigl(\Grp\bigl(\Alg(\fL)\bigr)\bigr)=\fL$, and therefore, $T$ is a collineation of
$\fL$, by \Cref{lem03}. This shows that $\Col(\fL)$ contains the normalizer of $\Grp\bigl(\Alg(\fL)\bigr)$.

(iv) Assume that $\Col(\fL)$ is the normalizer of $\Grp\bigl(\Alg(\fL)\bigr)$. Since $\Lat\Alg(\fL)$ is reflexive it
follows that $\Col\bigl(\Lat\Alg(\fL)\bigr)$ is the normalizer of $\Grp\bigl(\Alg\Lat\Alg(\fL)\bigr)=
\Grp\bigl(\Alg(\fL)\bigr)$, by (ii). Hence, by the assumption, we conclude $\Col\bigl(\Lat\Alg(\fL)\bigr)=\Col(\fL)$.
\end{proof}

\begin{corollary} \label{cor06}
An invertible operator $S\in \pB(\eX)$ is a collineation of a reflexive algebra $\pA$ if and only if $\Theta_S$ is a
spatial automorphism of $\pA$.
\end{corollary}

\begin{proof}
Let $S\in \pB(\eX)$ be invertible. Assume that $S\in\Col(\pA)$.
Since $\pA$ is reflexive, we have $\pA=\Alg\Lat(\pA)$. Hence, by the proof of \Cref{theo01}~(ii),
$S^{-1}\pA S=\pA$, that is $\Theta_S$ is a spatial automorphism of $\pA$. On the other hand,
if $S^{-1}\pA S=\pA$, then $AS\eM\subseteq S\eM$ and $AS^{-1}\eM\subseteq S^{-1}\eM$ for every
$A\in \pA$ and every $\eM\in \Lat(\pA)$. Hence, $S\in\Col(\pA)$, by \Cref{lem03}.
\end{proof}

The following theorem which characterizes collineations of reflexive algebras can be seen as an 
infinite-dimensional generalization of \cite[Theorem 2.5]{BK}.  

\begin{theorem} \label{theo02}
Let $\pA\subseteq \pB(\eX)$ be a reflexive algebra. For an invertible operator $S\in \pB(\eX)$, the following
assertions are equivalent.
\begin{enumerate}
\item[(a)] $S\in \Col(\pA)$;
\item[(b)] $S\pA_x\in \Lat(\pA)$ and $S^{-1}\pA_x\in \Lat(\pA)$ for all $x\in \eX$;
\item[(c)] $S\pA S^{-1}= \pA$;
\item[(d)] $\Lat(S\pA S^{-1})=\Lat(\pA)$;
\item[(e)] $S\pA_x=\pA_{Sx}$ for all $x\in \eX$.
\end{enumerate}
\end{theorem}

\begin{proof}
By \Cref{cor06}, (a) and (c) are equivalent, and, by \Cref{lem06}, (c) implies (e).
It is almost obvious that (e) implies (b). Indeed, for an arbitrary $y\in \eX$ replace $x$ with $S^{-1}y$ in 
the equality $S\pA_x=\pA_{Sx}$ to get $S^{-1}\pA_y=\pA_{S^{-1}y}$.
To prove that (b) implies (a), note first that cyclic subspaces generate $\Lat \pA$. Hence, subspaces $S\eM$ and
$S^{-1}\eM$ are in $\Lat(\pA)$ for every $\eM\in \Lat(\pA)$, so that $S\in\Col(\pA)$ by \Cref{lem03}.
It is clear that (c) implies (d). On the other hand, if (d) holds true, then $S\pA S^{-1}=
\Alg\Lat(S\pA S^{-1})=\Alg\Lat(\pA)=\pA$ since $S\pA S^{-1}$ is reflexive whenever $\pA$ is reflexive.
\end{proof}

For a subalgebra $\pA\subseteq \pB(\eX)$ and its commutant $\pA'$, it is easily seen that
$\Grp(\pA)\Grp(\pA')$ $=\{ AB;\, A\in \Grp(\pA),\, B\in \Grp(\pA')\}$ is a group generated by 
$\Grp(\pA)$ and $\Grp(\pA')$.

\begin{theorem} \label{theo05}
If $\pA\subseteq\pB(\eX)$ is a reflexive algebra, then $\Grp(\pA)\Grp(\pA')\subseteq \Col(\pA)\subseteq \Col(\pA')$.
Moreover, if $\pA$ and $\pA'$ are reflexive, then $\Grp(\pA)\Grp(\pA')$ is a normal subgroup of $\Col(\pA)$.
\end{theorem}

\begin{proof}
Let $AB\in \Grp(\pA)\Grp(\pA')$, where $A\in\Grp(\pA)$ and $B\in\Grp(\pA')$. By \Cref{lem03},
$A\in\Col(\pA)$. To see that $B\in\Col(\pA)$, choose an arbitrary $\eM\in \Lat(\pA)$.
Then for every $T\in \pA$ we have $TB\eM=BT\eM\subseteq B\eM$, that is, $B\eM\in \Lat(\pA)$.
Since $B^{-1}\in \Grp(\pA')$, we also have $B^{-1}\eM\in \Lat(\pA)$, so that $B\in \Col(\pA)$. Hence, 
$AB\in \Col(\pA)$, by \Cref{lem03}.

Suppose that $S\in \Col(\pA)$. We will show that $S^{-1}\pA' S=\pA'$. To this end, let $T\in\pA'$ be arbitrary.
Since $S^{-1}\pA S=\pA$, by \Cref{theo02}, for each $A\in\pA$ there exists $B\in\pA$ such that $A=S^{-1}BS$.
It follows that 
$$S^{-1}TSA=S^{-1}TBS=S^{-1}BTS=AS^{-1}TS.$$ 
This shows that $S^{-1}\pA' S\subseteq\pA'$. Since $S^{-1}\in\Col(\pA)$, we also have 
$S\pA' S^{-1}\subseteq\pA'$, and consequently $S^{-1}\pA' S=\pA'$. By \Cref{lem06}, 
$S^{-1}\bigl(\Alg\Lat(\pA')\bigr)S=\Alg\Lat(\pA')$, and therefore, $S\in \Col\bigl(\Alg\Lat(\pA')\bigr)$, 
by \Cref{theo02}. To finish the proof, notice that $\pA'$ and $\Alg\Lat(\pA')$ have the same lattice of 
invariant subspaces which implies $S\in\Col(\pA')$.

Assume now that $\pA$ and $\pA'$ are reflexive algebras. By \Cref{theo01}, $\Grp(\pA)$ and $\Grp(\pA')$ are 
normal subgroups of $\Col(\pA)$ and $\Col(\pA')$, respectively. Since $\Col(\pA)\subseteq\Col(\pA')$, 
for every $C\in\Col(\pA)$ we have $C\Grp(\pA)C^{-1}=\Grp(\pA)$ and $C\Grp(\pA')C^{-1}=\Grp(\pA')$.
Hence, $C\Grp(\pA)\Grp(\pA')C^{-1}=C\Grp(\pA)C^{-1} C\Grp(\pA')C^{-1}=\Grp(\pA)\Grp(\pA')$ which 
concludes the proof. 
\end{proof}

\begin{corollary} \label{cor01}
Let $\pT\subseteq\pB(\eX)$ be a non-empty set of operators.
\begin{itemize}
\item[(i)]  For $\pA=\Alg\Lat(\pT)$, we have $\Col(\pA)\subseteq \Col(\pA\cap\pA')$.
\item[(ii)] If $\pT'$ and $\pT''$ are reflexive algebras, then $\Col(\pT')=\Col(\pT'')$.
\end{itemize}
\end{corollary}

\begin{proof}
(i) Assume that $S\in \Col(\pA)$ and let $A\in \pA\cap\pA'$ be arbitrary. By \Cref{theo02}, $SAS^{-1}\in \pA$
and therefore $SAS^{-1}\in \pA'$, by the same reasoning as in the proof of \Cref{theo05}. Hence
$S(\pA\cap\pA')S^{-1}\subseteq \pA\cap\pA'$. Similarly, $S^{-1}(\pA\cap\pA')S\subseteq \pA\cap\pA'$
and therefore $S(\pA\cap\pA')S^{-1}=\pA\cap\pA'$. It follows, by \Cref{lem06},
$S\Alg\Lat\bigl(\pA\cap\pA'\bigr)S^{-1}= \Alg\Lat\bigl(\pA\cap\pA'\bigr)$
which means, by \Cref{theo02}, that $S\in \Col\bigl(\pA\cap\pA'\bigr)$.

(ii) Since $\pT'$ and $\pT''$ are reflexive algebras we have $\Col(\pT')\subseteq\Col(\pT'')\subseteq\Col(\pT''')$.
We conclude $\Col(\pT')=\Col(\pT'')$ as $\pT'''=\pT'$.
\end{proof}

Let $\eH$ be a complex Hilbert space. Recall that a von Neumann algebra is a self-adjoint \sot-closed
subalgebra $\pA\subseteq \pB(\eH)$ that contains $I$. By the von Neumann bicommutant theorem,
a self-adjoint subalgebra $\pA\subseteq \pB(\eH)$ is a von Neumann algebra if and only if $\pA=\pA''$. It follows
that $\pA'$ is a von Neumann algebra. It is well-known that every von Neumann algebra is reflexive.

\begin{corollary} \label{cor04}
If $\pA$ is a von Neumann algebra, then $\Col(\pA)=\Col(\pA')$.
\end{corollary}

\begin{proof}
Since $\pA'$ and $\pA''$ are reflexive algebras, by \Cref{cor01} we have $\Col(\pA')=\Col(\pA'')$, so that 
$\Col(\pA')=\Col(\pA)$ as $\pA$ is a von Neumann algebra.
\end{proof}

\subsection{A complement of $\Grp\bigl(\Alg(\fL)\bigr)$ in $\Col(\fL)$} \label{Sec033}

Let $\fL\subseteq \fC(\eX)$ be a subspace lattice. By \Cref{theo01}, $\Grp\bigl(\Alg(\fL)\bigr)$ is a
normal subgroup of $\Col(\fL)$. It is an interesting question whether $\Grp\bigl(\Alg(\fL)\bigr)$ is a complemented
subgroup of $\Col(\fL)$, that is, whether there exists a subgroup $\pO(\fL)\subseteq\Col(\fL)$ called a complement 
of $\Grp\bigl(\Alg(\fL)\bigr)$ in $\Col(\fL)$ such that $\Col(\fL)=\Grp\bigl(\Alg(\fL)\bigr)\pO(\fL)$ and 
$\Grp\bigl(\Alg(\fL)\bigr)\cap\pO(\fL)=\{ I\}$.
If this happens, then $\Col(\fL)$ is an inner semidirect product of $\Grp\bigl(\Alg(\fL)\bigr)$ and $\pO(\fL)$,
which we write as $\Col(\fL)=\Grp\bigl(\Alg(\fL)\bigr)\rtimes\pO(\fL)$. Moreover, since  $\Grp\bigl(\Alg(\fL)\bigr)$ 
is a normal subgroup of $\Col(\fL)$, we have 
$\Col(\fL)=\Grp\bigl(\Alg(\fL)\bigr)\pO(\fL)=\pO(\fL)\Grp\bigl(\Alg(\fL)\bigr)$. Indeed, if $S=AT$, with 
$A\in \Grp\bigl(\Alg(\fL)\bigr)$ and $T\in \pO(\fL)$, then $S=TB$, where $B=T^{-1}AT\in \Grp\bigl(\Alg(\fL)\bigr)$.

If a complement of $\Grp\bigl(\Alg(\fL)\bigr)$ in $\Col(\fL)$ exists, it does not need to be unique. 
Indeed, if $\pO(\fL)$ is a complement of $\Grp\bigl(\Alg(\fL)\bigr)$ and $\pO(\fL)$ is not a normal subgroup of 
$\Col(\fL)$, then there exists $C\in \Col(\fL)$ such that $C\pO(\fL)C^{-1}\ne \pO(\fL)$. It is not hard to see that 
$C \pO(\fL) C^{-1}$ is a complement of $\Grp(\fL)$, as well. Although complements of $\Grp(\fL)$ 
in $\Col(\fL)$ are possibly not unique, they are all isomorphic to the quotient group $\Col(\fL)/\Grp(\fL)$ via 
the natural group isomorphism $\pO(\fL)\to \Col(\fL)/\Grp(\fL)$ given by $T\mapsto T\Grp(\fL)$. 
Of special interest is the question of whether there exists a complement consisting of isometric operators.
Notice that this trivially holds if $\Col(\fL)=\Grp\bigl(\Alg(\fL)\bigr)$ in which case we have $\pO(\fL)=\{ I\}$.
The following simple proposition, whose proof is omitted, gives a necessary condition for the existence of a 
non-trivial complement of $\Grp\bigl(\Alg(\fL)\bigr)$ in $\Col(\fL)$.

\begin{proposition} \label{prop05}
Let $\fL\subseteq\fC(\eX)$ be a lattice of subspaces. The mapping $\phi\colon\Col(\fL)\to\Aut(\fL)$ given
by $\phi(S)=\Phi_S|_{\fL}$ is a group homomorphism  with kernel $\Grp\bigl(\Alg(\fL)\bigr)$. If $\Aut(\fL)$ is trivial,
then $\Col(\fL)=\Grp\bigl(\Alg(\fL)\bigr)$.\qed
\end{proposition}

\noindent
\begin{minipage}[c]{0.50\linewidth}
As we have seen, the triviality of $\Aut(\fL)$ implies $\Col(\fL)=\Grp\bigl(\Alg(\fL)\bigr)$. The converse does not hold.
For instance, let $\eM$ and $\eN$ be non-isomorphic subspaces of $\eX$ such that $\eM\wedge\eN=\{0\}$
and $\eM\vee\eN=\eX$. Then $\fM =\{\{0\},\eM,\eN ,\eX\}$ is a subspace lattice isomorphic to the diamond, see Figure 1,
whose group of automorphisms has two elements. However, $\Col(\fM)=\Grp\bigl(\Alg(\fM)\bigr)$.
\end{minipage}
\begin{minipage}[c]{0.5\linewidth}
\begin{center}
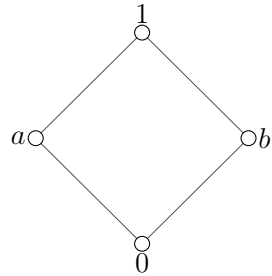

\begin{tikzpicture}[scale=1.4]
\draw[gray!70!black] (1cm,0cm)--(0cm,1cm)--(-1cm,0cm)--(0cm,-1cm)--cycle;
\fill[white] (1cm,0 cm) circle (2pt);
\fill[white] (-1cm,0 cm) circle (2pt);
\fill[white] (0cm,1 cm) circle (2pt);
\fill[white] (0cm,-1 cm) circle (2pt);
\draw (0cm,-1cm) circle (2pt) node[below, scale=0.9] {$0$};
\draw (-1cm,0cm) circle (2pt) node[left, scale=0.9] {$a$};
\draw (1cm,0cm) circle (2pt) node[right, scale=0.9] {$b$};
\draw (0cm,1cm) circle (2pt) node[above, scale=0.9] {$1$};
\end{tikzpicture}
\captionof{figure}{The diamond.} \label{fig01d}
\end{center}
\end{minipage}\smallskip

\noindent
Indeed, since $\eM$ and $\eN$ are not isomorphic it is impossible for an invertible operator $S\in \pB(\eX)$
to fulfill conditions $S\eM=\eN$ and $S\eN=\eM$.

A totally ordered lattice $C$ with at least one non-trivial element will be called a chain.
A lattice $M$ consisting of $n\geq 2$ chains
$C_i=\{ 0\neq a_{1}^{(i)}\leq a_{2}^{(i)}\leq\cdots\leq a_{q_i}^{(i)}\leq 1\}$ $(1\leq i\leq n)$
is a finite multi-chain if non-trivial elements $a_{k}^{(i)}\in C_i$ and $a_{l}^{(j)}\in C_j$ are incomparable
whenever $i\ne j$ (see \Cref{fig00}).

\noindent
\begin{minipage}[c]{0.50\linewidth}
\begin{center}
\begin{tikzpicture}[scale=1.2]
\draw[gray!70!black] (3cm,0cm)--(2.7cm,0.5cm);
\draw[gray!70!black] (2.7cm,0.5cm)--(2.7cm,1cm);
\draw[dotted,gray!70!black] (2.7cm,1cm)--(2.7cm,1.5cm);
\draw[gray!70!black] (2.7cm,1.5cm)--(3cm,2cm);
\fill[white] (2.7cm,0.5cm) circle (2pt);
\draw (2.7cm,0.5cm) circle (2pt) node[right, scale=0.9] {$a_{1}^{(2)}$};
\fill[white] (2.7cm,1cm) circle (2pt);
\draw (2.7cm,1cm) circle (2pt) node[right, scale=0.9] {$a_{2}^{(2)}$};
\fill[white] (2.7cm,1.5cm) circle (2pt);
\draw (2.7cm,1.5cm) circle (2pt) node[right, scale=0.9] {$a_{q_2}^{(2)}$};
\draw[gray!70!black] (3cm,0cm)--(2.2cm,0.5cm);
\draw[gray!70!black] (2.2cm,0.5cm)--(2.2cm,1cm);
\draw[dotted,gray!70!black] (2.2cm,1cm)--(2.2cm,1.5cm);
\draw[gray!70!black] (2.2cm,1.5cm)--(3cm,2cm);
\fill[white] (2.2cm,0.5cm) circle (2pt);
\draw (2.2cm,0.5cm) circle (2pt) node[left, scale=0.9] {$a_{1}^{(1)}$};
\fill[white] (2.2cm,1cm) circle (2pt);
\draw (2.2cm,1cm) circle (2pt) node[left, scale=0.9] {$a_{2}^{(1)}$};
\fill[white] (2.2cm,1.5cm) circle (2pt);
\draw (2.2cm,1.5cm) circle (2pt) node[left, scale=0.9] {$a_{q_1}^{(1)}$};
\draw[gray!70!black] (3cm,0cm)--(4cm,0.5cm);
\draw[gray!70!black] (4cm,0.5cm)--(4cm,1cm);
\draw[dotted,gray!70!black] (4cm,1cm)--(4cm,1.5cm);
\draw[gray!70!black] (4cm,1.5cm)--(3cm,2cm);
\fill[white] (4cm,0.5cm) circle (2pt);
\draw (4cm,0.5cm) circle (2pt) node[right, scale=0.9] {$a_{1}^{(n)}$};
\fill[white] (4cm,1cm) circle (2pt);
\draw (4cm,1cm) circle (2pt) node[right, scale=0.9] {$a_{2}^{(i)}$};
\fill[white] (4cm,1.5cm) circle (2pt);
\draw (4cm,1.5cm) circle (2pt) node[right, scale=0.9] {$a_{q_n}^{(n)}$};
\draw[dotted,gray!70!black] (3cm,0cm)--(3.4cm,0.5cm)--(3.4cm,1.5cm)--(3cm,2cm);
\draw[dotted,gray!70!black] (3cm,0cm)--(3.6cm,0.5cm)--(3.6cm,1.5cm)--(3cm,2cm);
\draw[dotted,gray!70!black] (3cm,0cm)--(3.8cm,0.5cm)--(3.8cm,1.5cm)--(3cm,2cm);
\fill[white] (3cm,0 cm) circle (2pt);
\draw (3cm,0 cm) circle (2pt) node[below, scale=0.9] {$0$};
\fill[white] (3cm,2cm) circle (2pt);
\draw (3cm,2cm) circle (2pt) node[above, scale=0.9] {$1$};
\end{tikzpicture}

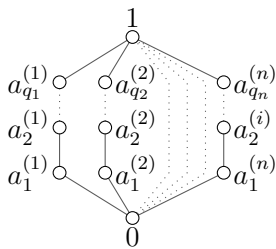
\captionof{figure}{A multi-chain.}   \label{fig00}
\end{center}
\end{minipage}
\begin{minipage}[c]{0.5\linewidth}
\begin{center}
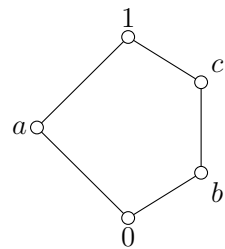

\begin{tikzpicture}[scale=1.2]
\draw (0cm,0cm)--(0.8cm,0.5cm);
\draw (0.8cm,0.5cm)--(0.8cm,1.5cm);
\draw (0.8cm,1.5cm)--(0cm,2cm);
\draw (0cm,0cm)--(-1cm,1cm);
\draw (-1cm,1cm)--(0cm,2cm);
\fill[white] (0cm,0 cm) circle (2pt);
\draw (0cm,0 cm) circle (2pt) node[below, scale=0.9] {$0$};
\fill[white] (-1cm,1 cm) circle (2pt);
\draw (-1cm,1cm) circle (2pt) node[left, scale=0.9] {$a$};
\fill[white] (0.8cm,0.5cm) circle (2pt);
\draw (0.8cm,0.5cm) circle (2pt) node[below right, scale=0.9] {$b$};
\fill[white] (0.8cm,1.5cm) circle (2pt);
\draw (0.8cm,1.5 cm) circle (2pt) node[above right, scale=0.9] {$c$};
\fill[white] (0cm,2cm) circle (2pt);
\draw (0cm,2 cm) circle (2pt) node[above, scale=0.9] {$1$};
\end{tikzpicture}
\captionof{figure}{The pentagon.}   \label{fig01p}
\end{center}
\end{minipage}\medskip

A lattice $L$ is realizable in a Banach space $\eX$ if there exists a lattice of subspaces
$\fL\subseteq \fC(\eX)$ and an isomorphism from $L$ to $\fL$. If this holds, then $\fL$ is a realization of $L$.
Not every lattice is realizable in every Banach space. However, it is easily seen that a chain with $q\in \mathbb N$
non-trivial elements can be realized as a subspace lattice in every Banach space of dimension at least $q+1$.
On the other hand, the pentagon, a lattice with five elements whose Hasse diagram is in \Cref{fig01p},
cannot be realized as a subspace lattice in a finite-dimensional Banach space, since such subspace lattices are
modular, whereas the pentagon is not (see \cite{Bra}).

A nest is well-ordered if every non-empty subset of it has a least subspace.
In the following theorem, we will show that well-ordered nests and some subspace lattices
that are realizations of finite multi-chains, have trivial groups of collineations.

\begin{theorem} \label{theo06}
If either
\begin{itemize}
\item[(i)] $\fM$ is a well-ordered nest or
\item[(ii)] $\fM$ is a realization of a finite multi-chain $M$ satisfying $q_i\ne q_j$ for all indices $i\ne j$,
\end{itemize}
then $\Col(\fM)=\Grp\bigl(\Alg(\fM)\bigr)$.
\end{theorem}

\begin{proof}
(i) Let $\fM$ be a well-ordered nest. We will prove that $\Aut(\fM)$ is trivial. To see this, let $\Phi\in\Aut(\fM)$
be arbitrary. Suppose that there exists a subspace $\eM\in \fM$ such that $\Phi(\eM)\ne \eM$. Hence,
$\fF=\{ \eK\in \fM;\; \Phi(\eK)\ne \eK\}$ is a non-empty subset of $\fM$.
By assumption, $\eM_0=\bigwedge\limits_{\eK\in \fF} \eK$ is the least subspace in $\fF$. Thus, either
$\Phi(\eM_0)\subsetneq \eM_0$ or $\Phi(\eM_0)\supsetneq \eM_0$. Assume that the former holds.
Since $\Phi(\eN)=\eN$ for every $\eN\in \fM$ such that $\eN\subsetneq \eM_0$, we have
$\Phi(\Phi(\eM_0))=\Phi(\eM_0)$ and therefore $\Phi(\eM_0)=\eM_0$ as $\Phi$ is injective. This contradicts
the fact that $\eM_0\in\fF$. If the latter holds, then $\Phi^{-1}(\eM_0)\subsetneq \eM_0$. As above, we can conclude that 
$\eM_0=\Phi(\Phi^{-1}(\eM_0))=\Phi^{-1}(\eM_0)$ yielding $\Phi(\eM_0)=\eM_0$. 

(ii) Let $\fM$ be a realization of a finite multi-chain $M$ satisfying $q_i\ne q_j$ for all indices $i\ne j$.
Hence, $\fM=\fC_1\cup\cdots\cup\fC_n$, where $\fC_i=\bigl\{ \{0\}\subsetneq \eM_{1}^{(i)}
\subsetneq \eM_{2}^{(i)}\subsetneq \cdots \subsetneq \eM_{q_i}^{(i)}\subsetneq \eX\bigr\}$ for $i=1,\ldots,n$.
There is no loss of generality if we assume that $q_1>\cdots>q_n$. We will show that $\Aut(\fM)$ is trivial.
Suppose that $\Phi\in\Aut(\fM)$. Then
$$\{0\}=\Phi(\{0\})\subsetneq \Phi(\eM_{1}^{(1)})\subsetneq \Phi(\eM_{2}^{(1)})\subsetneq\cdots
\subsetneq \Phi(\eM_{q_1}^{(1)})\subsetneq \Phi(\eX)=\eX$$ 
is a nest of length $q_1+2$ in $\fM$.
Since, $\fC_1$ is the only nest of length $q_1+2$ in $\fM$ we see that $\Phi$ preserves $\fC_1$, that is,
$\Phi(\eM_{j}^{(1)})\in\fC_1$ for every $j=1,\ldots,q_1$. However, as $\fC_1$ is a finite nest, it is necessarily 
well-ordered. Hence, by (i) we have $\Phi(\eM_{j}^{(1)})=\eM_{j}^{(1)}$ for every $j=1,\ldots,q_1$.
As the proof for the other chains follows in the same way, we omit it.
\end{proof}

\subsection{Adjoints of collineations} \label{Sec034}

The dual of a lattice $L$ with ordering $\leq_{_L}$ is the lattice $L^*$ with ordering
$\leq_{_{L^*}}$, where $L^*$ as a set is equal to $L$ and $x\leq_{_{L^*}} y$ if and only if $y\leq_{_{L}}x$,
for all $x, y\in L$. The ``identity mapping'' $\iota\colon L\to L^*$ is an anti-isomorphism of lattices.
The dual lattice of $L^*$ is clearly $L$. It is a natural question whether the dual $L^*$ of a lattice $L$, which
is realizable in a Banach space $\eX$ is realizable in its dual Banach space $\eX^*$. 

The pairing between
$\eX$ and $\eX^*$ is given by $\langle x,\xi\rangle=\xi(x)$ for all $x\in\eX$ and $\xi\in \eX^*$.
For $T\in\pB(\eX)$, let $T^*\in\pB(\eX^*)$ be its adjoint operator. If $\eU\subseteq \eX$ is a non-empty set, then let
$\eU^\perp=\{ \xi\in \eX^*;\; \langle x,\xi\rangle=0,\; \text{for every}\; x\in \eU\}$. Similarly, for a non-empty
set $\eV\subseteq \eX^*$, let $\eV_\perp=\{ x\in \eX;\; \langle x,\xi\rangle=0,\; \text{for every}\; \xi\in \eV\}$.
It is easily seen that $\eU^\perp$ is a subspace of $\eX^*$ which is closed in the weak$^*$ topology and
$\eV_\perp$ is a subspace of $\eX$.
For a non-empty family $\fF\subseteq \fC(\eX)$, let $\fF^\perp=\{ \eM^\perp\subseteq \eX^*;\; \eM\in \fF\}$.
Similarly, for $\emptyset\ne \fG\subseteq \fC(\eX^*)$, let $\fG_\perp=\{ \eN_\perp\subseteq \eX;\; \eN\in \fG\}$.

If $\fL$ is a realization of a lattice $L$ in a reflexive Banach space $\eX$, then $\fL^\perp$ is a realization
of the dual lattice $L^*$ in the dual space $\eX^*$. However, if $\eX$ is a non-reflexive Banach space, then
$\fL^\perp$ does need to be a subspace lattice in $\eX^*$.

\begin{example} \label{ex02}
Suppose that $\eX$ is not reflexive and $\eW\subseteq \eX^*$ is a subspace which is not closed in the
weak$^*$ topology (for instance, $\eW$ can be the kernel of a functional in $\eX^{**}\setminus\eX$). Let
$\{\eN_j;\; j\in J\}$ be a family of subspaces of $\eW$ which are closed in the weak$^*$ topology such that
$\bigvee_{j\in J}\eN_j=\eW$ (for instance, $\{\eN_j;\; j\in J\}$ can be the family of all one-dimensional subspaces of $\eW$). 
Let $\eM_j=(\eN_j)_\perp$ $(j\in J)$ and let $\fL$ be the subspace lattice generated by $\{ \eM_j;\; j\in J\}$.
It is clear that $\{ (\eM_j)^\perp;\; j\in J\}\subseteq \fL^\perp$ but $\eW=\bigvee_{j\in J}(\eM_j)^\perp$ is not in
$\fL^\perp$ since this family contains only subspaces which are closed in the weak$^*$ topology, yet $\eW$ is not
closed in the weak$^*$ topology.\qed
\end{example}

If $\pD\subseteq \pB(\eX)$ is a non-empty set, let $\pD^*=\{ T^*\in \pB(\eX^*);\; T\in\pD\}$. For instance,
if $\fF\subseteq \fC(\eX)$ is a non-empty family, then $\Col(\fF)^*=\{ S^*\in\pB(\eX^*);\; S\in \Col(\fF)\}$.
It is clear that $\Col(\fF)^*$ is a subgroup of the group of invertible operators on $\eX^*$.

\begin{proposition} \label{prop02}
Let $\fF\subseteq \fC(\eX)$ be a non-empty family. Then
\begin{itemize}
\item[(i)] $\Col(\fF)^*$ is a subgroup of $\Col(\fF^\perp)$ and
\item[(ii)] $\Col(\fF^\perp)=\Col(\fF)^*$ if $\eX$ is reflexive.
\end{itemize}
\end{proposition}

\begin{proof}
(i) Let $S\in \Col(\fF)$ and $\eM\in \fF$ be arbitrary. If $\xi\in (S^{-1}\eM)^\perp$, then
$\langle x, (S^{-1})^*\xi\rangle=\langle S^{-1}x,\xi\rangle=0$ for all $x\in \eM$. Hence,
$(S^{-1})^*\xi\in \eM^\perp$ and therefore $\xi\in S^*\eM^\perp$. This shows that
$(S^{-1}\eM)^\perp\subseteq S^*\eM^\perp$. To prove the opposite inclusion, let $\eta\in \eM^\perp$.
For every $y\in S^{-1}\eM$, there exists $x\in \eM$ such that $y=S^{-1}x$. Hence,
$\langle y,S^*\eta\rangle=\langle S S^{-1}x,\eta\rangle=0$. We have proved that
$ S^*\eM^\perp=(S^{-1}\eM)^\perp\in\fF^\perp$. Similarly, $ (S^*)^{-1}\eM^\perp=(S\eM)^\perp\in\fF^\perp$.
Thus, $S^*\in \Col(\fF^\perp)$, by \Cref{lem03}.

(ii) If $\eX$ is reflexive, then we may identify $\eX$ and $\eX^{**}$. In this sense, we have $T^{**}=T$ for every
operator $T\in\pB(\eX)$ and $(\fF^\perp)^\perp=\fF$. Using (i) we see
that $\Col(\fF)= \Col(\fF)^{**}\subseteq\Col(\fF^\perp)^*\subseteq \Col((\fF^{\perp})^{\perp})=\Col(\fF)$
which gives $\Col(\fF^\perp)=\Col(\fF)^*$.
\end{proof}

A subspace $\eM\subseteq\eX$ is invariant for an operator $T\in \pB(\eX)$ if and only if the subspace
$\eM^\perp\subseteq\eX^*$ is invariant for the adjoint operator $T^*$. Hence, if $\eX$ is a reflexive Banach space
and $\fL\subseteq \fC(\eX)$ is a subspace lattice, then $\Alg(\fL^\perp)=\Alg(\fL)^*$ and therefore
$\Grp\bigl(\Alg(\fL^\perp)\bigr)=\Grp\bigl(\Alg(\fL)\bigr)^*$. Suppose that
$\Col(\fL)=\Grp\bigl(\Alg(\fL)\bigr)\rtimes\pO(\fL)$. By \Cref{prop02}~(ii), every operator in
$\Col(\fF^\perp)$ is of the form $S^*$ for some operator $S\in\Col(\fL)$. Let $A\in \Grp\bigl(\Alg(\fL)\bigr)$
and $V\in\pO(\fL)$ be such that $S=A V$. Then $S^*=(V^{-1}AV)^* V^*$, where
$(V^{-1}AV)^*\in \Grp\bigl(\Alg(\fL)\bigr)^*$ (recall that $\Grp\bigl(\Alg(\fL)$ is a normal subgroup of $\Col(\fL)$)
and $V^*\in \pO(\fL)^*$. Since $\Grp\bigl(\Alg(\fL)\bigr)^*\cap \pO(\fL)^*$ contains only the identity operator
on $\eX^*$ we conclude that the following holds.

\begin{corollary} \label{cor07}
Let $\eX$ be a reflexive Banach space. If $\fL\subseteq \fC(\eX)$ is a subspace lattice such that
$\Col(\fL)=\Grp\bigl(\Alg(\fL)\bigr)\rtimes\pO(\fL)$, then
$\Col(\fL^\perp)=\Grp\bigl(\Alg(\fL^\perp)\bigr)\rtimes\pO(\fL^\perp)$, where
$\Grp\bigl(\Alg(\fL^\perp)\bigr)=\Grp\bigl(\Alg(\pL)\bigr)^*$ and $\pO(\fL^\perp)=\pO(\fL)^*$.\qed
\end{corollary}

\section{Conjugate collineations} \label{Sec04}
\setcounter{theorem}{0}
Let $\eX$ be a complex Banach space. Recall that a continuous mapping $\bar{T}\colon\eX\to\eX$ is a
conjugate linear operator if $\bar{T}(\alpha x+\beta y)=\overline{\alpha}\bar{T}x+\overline{\beta}\bar{T}y$
for arbitrary $x, y\in \eX$, $\alpha, \beta\in\bC$. The norm of $\bar{T}$ is
$\|\bar{T}\|=\sup\{\|\bar{T}x\|;\; x\in\eX,\; \|x\|\leq 1\}$. It is not hard to see that $\overline{\pB}(\eX)$,
the set of all conjugate linear operators on $\eX$, is a complex Banach space.
We use a bar to stress that operators in $\overline{\pB}(\eX)$ are conjugate linear. The same convention is used
for subsets of $\overline{\pB}(\eX)$. The Banach space
$\overline{\pB}(\eX)$ is not an algebra because the product of two conjugate linear operators is not conjugate linear.
The definition of an invertible conjugate linear operator $\bar{T}$ is obvious: there must exists a conjugate linear
operator $\bar{T}^{-1}$ such that $\bar{T}^{-1}\bar{T}=\bar{T}\bar{T}^{-1}=I$. As we shall see later, there
exist complex Banach spaces without invertible conjugate linear operators.

Every invertible conjugate linear operator $\bar{S}$ on $\eX$ induces an automorphism
$\Phi_{\bar{S}}$ of $\fC(\eX)$. Let $\fL\subseteq \fC(\eX)$ be a subspace lattice. Analogous to \Cref{def01},
we will say that  an invertible conjugate linear operator $\bar{S}$ is a {\em conjugate collineation} of $\fL$ if
a subspace $\eM$ belongs to $\fL$ if and only if $\bar{S}\eM$ belongs to $\fL$. The set of all conjugate collineations
of $\fL$ will be denoted by $\overline{\Col}(\fL)$. In the case when $\fL=\Lat(\pA)$, for an algebra $\pA$, we
will write $\overline{\Col}(\pA)$ instead of $\overline{\Col}(\Lat(\pA))$.

For a non-empty family $\fF\subseteq \fC(\eX)$, let $\overline{\Alg}(\fF)=\{\bar{T}\in \overline{\pB}(\eX);\;
\bar{T}\eM\subseteq \eM,\; \text{for all}\; \eM\in \fF\}$.
It is clear that $\overline{\Alg}(\fF)$ is a \sot-closed subspace of $\overline{\pB}(\eX)$. Moreover, it is not hard to
see that it is a Banach $\Alg(\fF)$-bimodule. For non-empty subsets $\overline{\pS}_1$ and $\overline{\pS}_2$
of $\overline{\pB}(\eX)$, let $\overline{\pS}_1\cdot\overline{\pS}_2$ denote the closed linear span of
$\{ \bar{S}_1\bar{S}_2;\; \bar{S}_1\in \overline{\pS}_1,\bar{S}_2\in \overline{\pS}_2\}$. In particular,
for a non-empty subset $\overline{\pS}\subseteq \overline{\pB}(\eX)$, let
$\overline{\pS}^{[2]}=\overline{\pS}\cdot \overline{\pS}$. It is clear that these are subspaces of $\pB(\eX)$.
It is easily seen that $\overline{\Alg}(\fF)^{[2]}\subseteq \Alg(\fF)$.

\begin{proposition} \label{prop07}
Let $\fL\subseteq \fC(\eX)$ be a subspace lattice.
\begin{itemize}
\item[(i)] For an invertible conjugate linear operator $\bar{S}$, the following are equivalent:
	\begin{itemize}
	\item[(a)] $\bar{S}\in \overline{\Col}(\fL)$;
	\item[(b)] if $\eM\in\fL$, then $\bar{S}\eM$ and $\bar{S}^{-1}\eM$ are in $\fL$;
	\item[(c)] $\bar{S}^{-1}$ is a collineation of $\fL$.
	\end{itemize}
\item[(ii)] If $\bar{S}\in \overline{\Alg}(\fL)$ is invertible and $\bar{S}^{-1}\in \overline{\Alg}(\fL)$,
then $\bar{S}\in \overline{\Col}(\fL)$.
\item[(iii)] If $\bar{S}$ is an invertible conjugate linear operator commuting with every operator from $\Alg(\fL)$, 
then $\bar{S}\in \overline{\Col}(\Alg(\fL))$.
\item[(iv)] $\overline{\Col}(\fL)\subseteq \overline{\Col}(\Alg(\fL))$.
\end{itemize}
\end{proposition}

\begin{proof}
(i) This is an analog of \Cref{lem03} and the proof is almost the same. (ii) follows from (i). To prove (iii),
let $\bar{S}$ be an invertible conjugate linear operator commuting with every operator from $\Alg(\fL)$. 
Then the same holds for its inverse $\bar{S}^{-1}$.
If $\eM\in \Lat\Alg(\fL)$, then $A\bar{S}\eM=\bar{S}A\eM\subseteq \bar{S}\eM$ for all $A\in \Alg(\fL)$. Hence,
$\bar{S}\eM\in \Lat\Alg(\fL)$. Similarly, $\bar{S}^{-1}\eM\in \Lat\Alg(\fL)$ for every $\eM\in \Lat\Alg(\fL)$.
By (i), $\bar{S}\in \overline{\Col}(\Lat\Alg(\fL))$. The proof of (iv) is a simple modification of the proof of
\Cref{theo01}~(ii).
\end{proof}

It follows from \Cref{prop07} that $\overline{\Col}(\fL)$ is not empty if there exists an invertible operator
$\bar{S}\in \overline{\Alg}(\fL)$ such that  $\bar{S}^{-1}\in \overline{\Alg}(\fL)$. Moreover, the
following proposition shows that in this case, conjugate collineations can be easily derived from collineations.

\begin{proposition} \label{prop10}
If there exists a conjugate collineation $\bar{S}\in \overline{\pB}(\eX)$ of a subspace lattice
$\fL\subseteq \fC(\eX)$, then $ \overline{\Col}(\fL)=\bar{S}\Col(\fL)=\Col(\fL)\bar{S}$ and
$\Col(\fL)=\bar{S}\, \overline{\Col}(\fL)=\overline{\Col}(\fL)\bar{S}$.
\end{proposition}

\begin{proof}
Suppose that there exists $\bar{S}$ in $\overline{\Col}(\fL)$. Then we have
$\bar{T}=\bar{S}\bar{S}^{-1}\bar{T}=\bar{T}\bar{S}^{-1}\bar{S}$ for every $\bar{T}\in \overline{\Col}(\fL)$.
Notice that $\bar{S}^{-1}\bar{T}$ and $\bar{T}\bar{S}^{-1}$ are invertible linear operators. Since $\bar{S}$
and $\bar{T}$ are in $\overline{\Col}(\fL)$ it is obvious that $\bar{S}^{-1}\bar{T}\eM\in \fL$ and
$\bar{T}^{-1}\bar{S}\eM\in \fL$ for all $\eM\in \fL$. It follows that
$\bar{S}^{-1}\bar{T}$ is a collineation of $\fL$. Hence, $\overline{\Col}(\fL)\subseteq \bar{S}\Col(\fL)$.
A similar reasoning gives $\overline{\Col}(\fL)\subseteq \Col(\fL)\bar{S}$.

To prove the opposite inclusions, let $\bar{S}\in \overline{\Col}(\fL)$ and $T\in \Col(\fL)$. Then
$\bar{S}T$ and $(\bar{S}T)^{-1}$ are invertible conjugate linear operators. Since $\bar{S}T\eM\in \fL$ and
$(\bar{S}T)^{-1}\eM=T^{-1}\bar{S}^{-1}\eM\in \fL$ for all $\eM\in \fL$, we conclude that
$\bar{S}\Col(\fL)\subseteq \overline{\Col}(\fL)$. Similarly, $\Col(\fL)\bar{S}\subseteq \overline{\Col}(\fL)$.
This proves the equalities $ \overline{\Col}(\fL)=\bar{S}\Col(\fL)=\Col(\fL)\bar{S}$.

Since $\bar{S}\in \overline{\Col}(\fL)$ if and only if $\bar{S}^{-1}\in \overline{\Col}(\fL)$, the first part of this
proof gives $ \overline{\Col}(\fL)=\bar{S}^{-1}\Col(\fL)=\Col(\fL)\bar{S}^{-1}$. It follows that
$\Col(\fL)=\bar{S}\, \overline{\Col}(\fL)=\overline{\Col}(\fL)\bar{S}$.
\end{proof}

There may be no conjugate linear conjugation of a subspace lattice. Even more is true. There exist Banach spaces
without invertible conjugate linear operators on them. To see this, we will use a very deep result proved by
Bourgain \cite{Bou}. Of course, classical Banach spaces have invertible conjugate linear operators. For instance,
the mapping $f(t)\mapsto \overline{f(t)}$ is a conjugate linear involution on several Banach spaces of functions.

The complex conjugate $\overline{\eX}$ of a complex Banach space $\eX$ is the Banach space with the same 
norm and addition as $\eX$, but with scalar multiplication defined by $(\alpha, x) \mapsto \bar{\alpha} x$ 
(see \cite{Bou, Kal}).
Banach spaces $\eX$ and $\overline{\eX}$ are isometrically isomorphic as real Banach spaces.
Let $J\colon \eX \to \overline{\eX}$ be the mapping given by $J x=x$ for all $x\in \eX$.
It is easily seen that $J$ is a conjugate linear bijective isometry. The same holds for its inverse
$J^{-1}\colon \overline{\eX}\to \eX$. If $T$ is a linear operator from $\eX$ to $\overline{\eX}$,
that is $T\in \pB(\eX,\overline{\eX})$, then
$J^{-1}T(\alpha x)=J^{-1}\bigl( \alpha\cdot Tx\bigr)=\overline{\alpha}Tx$ for all $\alpha\in \bC$, $x\in \eX$.
Hence, $J^{-1}T$ is a conjugate linear operator on $\eX$. Similarly, if $\bar{S}\in \overline{\pB}(\eX)$,
then $J\bar{S}\in\pB(\eX,\overline{\eX})$.

\begin{proposition} \label{prop11}
Mapping $\iota\colon \overline{\pB}(\eX)\to \pB(\eX,\overline{\eX})$ which is given by $\iota(\bar{S})=J\bar{S}$
is an isometric conjugate linear bijection. A conjugate linear operator $\bar{S}\in \overline{\pB}(\eX)$ is
invertible if and only if $\iota(\bar{S})$ is an isomorphism of complex Banach spaces $\eX$ and $\overline{\eX}$.
\end{proposition}

\begin{proof}
It is straightforward to check that $\iota$ is an isometric conjugate linear bijection.
To verify the second statement, assume that $\bar{S}\in \overline{\pB}(\eX)$ is invertible.
Then $\iota(\bar{S})=J\bar{S}$ is a bijection and therefore $\iota(\bar{S})$ is an isomorphism.
On the other hand, if $\iota(\bar{S})=J\bar{S}$ is an isomorphism, then $\bar{S}=J^{-1}\iota(\bar{S})$ is
invertible since $J^{-1}$ and $\iota(\bar{S})$ are bijective.
\end{proof}

\begin{corollary} \label{cor03}
There exists a complex Banach space $\eX_B$ without invertible conjugate linear operators.
\end{corollary}

\begin{proof}
Bourgain \cite{Bou} showed that there exists a Banach space $\eX_B$ such that $\eX_B$ and its complex conjugate
$\overline{\eX_{B}}$ are not isomorphic as complex Banach spaces (for a more elementary example see \cite{Kal}).
By \Cref{prop11}, it follows that there is no invertible conjugate linear operator on $\eX_B$.
\end{proof}

The following example shows that even in the case when there exist invertible conjugate linear operators on $\eX$
it is possible that $\overline{\Col}(\fL)$ is empty for a subspace lattice $\fL\subseteq \fC(\eX)$.

\begin{example} \label{ex01}
Let $\eX_B$ be a complex Banach space without invertible conjugate linear operators and let
$\eX=\eX_B\oplus\overline{\eX_B}$. Recall that $J\colon \eX_B\to\overline{\eX_B}$ is given by $Jx=x$ for all
$x\in \eX_B$. It is not hard to see that $\bar{V}\colon \eX\to \eX$, which is given by
$\bar{V}(x\oplus y)=J^{-1}y\oplus Jx$, is an invertible conjugate linear operator on $\eX$.
It follows that $\bar{V}T$ is an invertible conjugate linear operator for every invertible operator $T\in\pB(\eX)$.
Let $\eM=\eX_B\oplus\{0\}$ and
$\fL=\bigl\{ \{0\},\eM,\eX\bigr\}\subseteq \fC(\eX)$. We claim that $\overline{\Col}(\fL)$ is empty. Indeed, if
there were $\bar{S}$ in $\overline{\Col}(\fL)$, then we would have $\bar{S}\eM= \eM$ and therefore
the restriction $\bar{S}|_{\eM}$ would be an invertible conjugate linear operator.
But this is impossible since there is no invertible conjugate linear operator on $\eM$.\qed
\end{example}

\section{Collineations of some subspace lattices} \label{Sec05}
\setcounter{theorem}{0}

\subsection{Medial lattices} \label{Sec051}
A multi-chain is a medial lattice if it consists of $n\geq 2$ chains and each chain has exactly one
non-trivial element. Hence, in a medial lattice, every non-trivial element $a$ has at least one complement, that is,
there exists a non-trivial element $a'$ in the lattice such that $a\wedge a'=0$ and $a\vee a'=1$.\medskip

\noindent
\begin{minipage}[c]{0.450\linewidth}
The simplest medial lattices are the diamond and the double triangle (see Figures \ref{fig01d} and \ref{fig01t}). 
Note that the group of automorphisms of a medial lattice with $n$ non-trivial elements is isomorphic to the
symmetric group $S_n$.
\end{minipage}
\begin{minipage}[c]{0.55\linewidth}
\begin{center}
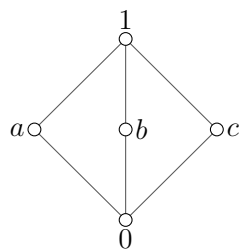

\begin{tikzpicture}[scale=1.2]
\draw[gray!70!black] (1cm,0cm)--(0cm,1cm)--(-1cm,0cm)--(0cm,-1cm)--cycle;
\draw[gray!70!black] (0cm,-1cm)--(0cm,1cm);
\fill[white] (0cm,0 cm) circle (2pt);
\fill[white] (1cm,0 cm) circle (2pt);
\fill[white] (-1cm,0 cm) circle (2pt);
\fill[white] (0cm,1 cm) circle (2pt);
\fill[white] (0cm,-1cm) circle (2pt);
\draw (0cm,-1cm) circle (2pt) node[below, scale=0.9] {$0$};
\draw (1cm,0cm) circle (2pt) node[right, scale=0.9] {$c$};
\draw (-1cm,0cm) circle (2pt) node[left, scale=0.9] {$a$};
\draw (0cm,0 cm) circle (2pt) node[right, scale=0.9] {$b$};
\draw (0cm,1cm) circle (2pt) node[above, scale=0.9] {$1$};
\end{tikzpicture}
\captionof{figure}{The double triangle.} \label{fig01t}
\end{center}
\end{minipage}\smallskip

It is not hard to see that the diamond can be realized in every Banach space of dimension at least $2$.
On the other hand, the double triangle can be realized in a Banach space $\eX$ if there exist isomorphic subspaces
$\eM,\eN\in \fC(\eX)$ that are algebraically complemented, that is, $\eM\cap\eN=\{0\}$ and $\eM+\eN=\eX$
(briefly, $\eM\oplus\eN=\eX$). In particular, the double triangle can be realized in a finite-dimensional Banach
space of dimension $2n$ $(n\in\bN)$ and in every infinite-dimensional Hilbert space (see \cite{Bra}).

\begin{proposition} \label{prop04}
Let $\eX$ be a complex Banach space and let $\eM,\eN\in \fC(\eX)$ be such that $\eX=\eM\oplus\eN$. Then $\eM$ 
and $\eN$ are isomorphic if and only if there exists a subspace $\eK\subseteq \eX$ distinct from both $\eM$ and 
$\eN$ such that $\eM\oplus\eK=\eX=\eN\oplus\eK$. In particular, $\eK$ is isomorphic to $\eM$ and to $\eN$.
\end{proposition}
 
\begin{proof}
Suppose that $\eM$ and $\eN$ are isomorphic and let $V\colon \eM\to\eN$ be an isomorphism. It is obvious that
$\eK:=\{ x+Vx;\; x\in\eM\}=\{ V^{-1}y+y;\; y\in \eN\}$
is a linear submanifold of $\eX$. Since $\eX=\eM\oplus \eN$ and both $\eM$ and $\eN$ are closed, the projections
onto $\eM$ and $\eN$ are continuous. Using the fact that $V$ is an isomorphism it is straightforward to verify that 
$\eK$ is closed. If $x\in \eM\cap \eK$, then there exists $x'\in \eM$ such that $x=x'+Vx'$. Hence, 
$x-x'=Vx'\in \eM\cap\eN=\{0\}$, that is, $x=0$.
If $z\in \eX$, then there exist unique $x\in\eM$ and $y\in \eN$ such that $z=x+y$. 
Since $x-V^{-1}y\in \eM$ and $V^{-1}y+y\in\eK$, from the identity $z=x-V^{-1}y+V^{-1}y+y$ we can conclude that
$\eX=\eM\oplus\eK$. Similarly, $\eX=\eN\oplus\eK$.

Assume that there exists a subspace $\eK\subseteq \eX$ such that $\eM\oplus\eK=\eX=\eN\oplus\eK$. 
Since $\eM$ and $\eN$ are complementary subspaces of a subspace $\eK$ in a Banach space $\eX$, 
by \cite[Corollary 3.2.16]{Meg}  they are isomorphic.
\end{proof}

Let $\eM, \eN$ be subspaces of $\eX$ such that $\eX=\eM\oplus\eN$. If $x\in\eM$ and $y\in\eN$, then
we will write $x+y$ instead of $x\oplus y$ since $\eM$ and $\eN$ are subspaces of $\eX$. It is well-known that
operators $A\in\pB(\eM)$, $B\in \pB(\eN,\eM)$, $C\in \pB(\eM,\eN)$, and $D\in \pB(\eN)$ induce an operator on 
$\eX=\eM\oplus\eN$ given by $x+y\mapsto (Ax+By)+(Cx+Dy)$ for all $x\in \eM$ and $y\in \eN$.
We will omit a simple proof of the following corollary.

\begin{corollary} \label{cor02}
If $\eK,\eM$ and $\eN$ are subspaces of $\eX$ such that $\eX=\eM\oplus\eN=\eM\oplus\eK=\eN\oplus\eK$
and $V\colon\eM\to\eN$ is an isomorphism such that $\eK=\{ x+Vx;\; x\in\eM\}$, then the operator
$W\colon x+y\mapsto V^{-1}y+Vx$, where $x\in\eM$ and $y\in\eN$ are arbitrary, has the
following properties: $W^2=I$, $W\eM=\eN$, $W\eN=\eM$, and $Wz=z$ for every $z\in \eK$, that is,
$\eK=\Ker(I-W)$. \qed
\end{corollary}

Let $\eM_1,\eM_2,\eM_3\in\fC(\eX)$ be such that $\eX=\eM_1\oplus\eM_2=\eM_1\oplus\eM_3=
\eM_2\oplus\eM_3$ and let $V_1\colon\eM_2\to\eM_3$, $V_2\colon\eM_1\to\eM_3$, and
$V_3\colon \eM_1\to\eM_2$ be isomorphisms such that
$ \eM_1=\{ x+V_1 x;\; x\in \eM_2\}=\{ V_{1}^{-1}y+y;\; y\in \eM_3\},$
$ \eM_2=\{ x+V_2 x;\; x\in \eM_1\}=\{ V_{2}^{-1}y+y;\; y\in \eM_3\},$ and
$ \eM_3=\{ x+V_3 x;\; x\in \eM_1\}=\{ V_{3}^{-1}y+y;\; y\in \eM_2\}.$
It is obvious that $\fD=\bigl\{ \{ 0\},\eM_1,\eM_2,\eX\bigr\}$ and
$\fT=\bigl\{ \{ 0\},\eM_1,\eM_2,\eM_3,\eX\bigr\}$ are realizations  of the diamond and the double
triangle, respectively. Define operators $W_1, W_2, W_3\in \pB(\eX)$ by
\begin{equation} \label{eq04}
\begin{split}
W_1(x_2+x_3)&=V_{1}^{-1}x_3+V_1x_2,\qquad (x_2\in \eM_2, x_3\in\eM_3),\\
W_2(x_1+x_3)&=V_{2}^{-1}x_3+V_2x_1,\qquad (x_1\in \eM_1, x_3\in\eM_3),\\
W_3(x_1+x_2)&=V_{3}^{-1}x_2+V_3x_1,\qquad (x_1\in \eM_1, x_2\in\eM_2).
\end{split}
\end{equation}
Although, in general, operators $W_1, W_2$, and $W_3$ are not isometric, by \Cref{cor02} they satisfy 
$W_{1}^{2}=W_{2}^{2}=W_{3}^{2}=I$. By the same corollary, $W_3$ is a collineation of $\fD$, and 
$W_1, W_2, W_3$ are collineations of $\fT$. Let us denote $\pO(\fD)=\{ I,W_3\}$ and 
$\pO(\fT)=\{I, W_1,W_2,W_3,W_1W_2,W_2W_1\}$.  It is obvious that $\pO(\fD)$ is a subgroup of $\Col(\fD)$
isomorphic to the symmetric group $S_2$.

\begin{lemma} \label{lem10}
$\pO(\fT)$ is a subgroup of $\Col(\fT)$ isomorphic to the symmetric group $S_3$.
\end{lemma}

\begin{proof}
Let $V_j$ and $W_j$ $(j=1,2,3)$ be as above. First, we will show that
$V_3=V_{1}^{-1} V_{2}$.
If $x_1\in \eM_1$, then there exists a unique vector $x_2\in \eM_2$ such that $x_1=x_2+V_1 x_2$.
Similarly, there exists a unique vector $y_{1}\in \eM_1$ such that $x_2=y_{1}+V_2 y_{1}$. Hence,
$x_2=x_1-V_1 x_2=y_1+V_2y_1$. Since $\eX=\eM_1\oplus \eM_3$ we conclude that $y_{1}=x_1$ and
$V_2 y_1=-V_1 x_2$, and therefore, $x_2=-V_{1}^{-1}V_2 x_1$. Since $y_1=x_1$ we have $x_2=x_1+V_2 x_1$. 
On the other hand, there exists a unique vector $z_{1}\in \eM_1$ such that
$V_2 x_1=z_{1}+V_3 z_{1}$. As above, since $\eX=\eM_1\oplus \eM_2$, from
$V_2x_1=-x_1+x_2=z_{1}+V_3 z_{1}$ we conclude that $z_{1}=-x_1$ and $x_2=V_3 z_{1}=-V_3x_1$.
This proves that $V_{1}^{-1}V_2 x_1=V_3x_1$ holds for an arbitrary vector $x_1\in \eM_1$.

Now we will show that
\begin{equation} \label{eq07}
W_1W_2=W_3 W_1=W_2W_3\qquad \text{and}\qquad W_1W_3=W_2W_1=W_3 W_2.
\end{equation}
We will check only the equality $W_1 W_2=W_3 W_1$ as other equalities can be proved similarly.
Let $y\in \eX$ be arbitrary. Then there exist unique vectors $x_1\in \eM_1$ and $x_2\in \eM_2$ such that
$y=x_1+x_2$. Hence, $W_1W_2 y=W_1(W_2 x_1+W_2 x_2)=W_1(W_2 x_1+x_2)=V_{1}^{-1}V_2 x_1+V_1 x_2$,
where the last equality holds because of \eqref{eq04}. On the other hand, $W_3 W_1 y=
W_3(W_1 x_1+W_1 x_2)=W_3(x_1+V_1 x_2)=V_3 x_1+V_1 x_2$.
Since $V_3=V_{1}^{-1}V_2$, by the first paragraph of this proof, we conclude that $W_1 W_2=W_3 W_1$.

Relations from \eqref{eq07} together with $W_{1}^{2}=W_{2}^{2}=W_{3}^{2}=I$ imply that $\pO(\fT)$ is
a subgroup of $\Col(\fT)$. Since $\pO(\fT)$ is a non-abelian group of order $6$ it is isomorphic to the symmetric
group $S_3$.
\end{proof}

\begin{theorem} \label{theo04}
Let $\eM_1,\eM_2,\eM_3\in\fC(\eX)$ be such that $\eX=\eM_1\oplus\eM_2=\eM_1\oplus\eM_3=
\eM_2\oplus\eM_3$ and let $\fD=\bigl\{ \{ 0\},\eM_1,\eM_2,\eX\bigr\}$ and
$\fT=\bigl\{ \{ 0\},\eM_1,\eM_2,\eM_3,\eX\bigr\}$ be realizations of the diamond and the double
triangle, respectively. Let $W_1, W_2, W_3\in \pB(\eX)$ be given by \eqref{eq04}. Then 
$$\Col(\fD)=\Grp\bigl(\Alg(\fD)\bigr)\rtimes \pO(\fD)\quad\text{and}\quad
\Col(\fT)=\Grp\bigl(\Alg(\fT)\bigr)\rtimes \pO(\fT),$$
where $\pO(\fD)=\{ I,W_3\}$ and $\pO(\fT)=\{I, W_1,W_2,W_3,W_1W_2,W_2W_1\}$.
\end{theorem}

\begin{proof}
It is clear that $\Grp\bigl(\Alg(\fT)\bigr)\cap \pO(\fT)=\{ I\}$. Let $S$ be an arbitrary collineation of $\fT$ and let
$\pi$ be a permutation of $\{ 1,2,3\}$ such that $S\eM_{\pi(j)}=\eM_j$ for $j=1,2,3$. If $W\in \pO(\fT)$
is the operator for which $W\eM_j=\eM_{\pi(j)}$, then $SW\eM_j=\eM_j$ for $j=1,2,3$. Hence,
$SW\in \Grp\bigl(\Alg(\fT)\bigr)$ and $S=SW W^{-1}$. This proves that
$\pO(\fT)=\{I, W_1,W_2,W_3,W_1W_2,W_2W_1\}$ is a complement of $\Grp\bigl(\Alg(\fT)\bigr)$ in $\Col(\fT)$.
The proof that $\pO(\fD)=\{ I,W_3\}$ is a complement of $\Grp\bigl(\Alg(\fD)\bigr)$ in $\Col(\fD)$ is similar.
\end{proof}

\subsection{Nests} \label{Sec052}
Recall that in this paper, a nest means a totally ordered subspace lattice in a Banach space $\eX$.
We refer the reader to the classical monograph \cite{Dav} about nests in Hilbert spaces.
Every nest is a reflexive subspace lattice (see \cite[Theorem 3.4]{Rin1} and \cite[Theorem 4.4]{Bra}).
By \Cref{theo06}, $\Col(\fN)=\Grp\bigl(\Alg(\fN)\bigr)$ if $\fN\subseteq\fC(\eX)$ is a well-ordered nest.
Now we will consider one of the simplest nests that is not well-ordered.

For $1\leq p< \infty$, let $\ell^p(\bZ)$ be the classical Banach space of two-sided complex sequences of the form 
$x=(\alpha_j)_{j\in\bZ}$ such that
$\| x\|_p=\bigl(\sum_{j\in\bZ} |\alpha_j|^p\bigr)^{1/p}<\infty$. Denote by $(e_j)_{j\in\bZ}$ the standard basis
of $\ell^p(\bZ)$. For each $k\in \bZ$, let $\eN_k=\bigvee\{ e_j;\; j\leq k\}$. Then
$\fN=\{ \eN_k;\; k\in \bZ\}\cup\bigl\{ \{0\},\ell^p(\bZ)\bigr\}$ is a nest.
A bilateral shift on $\ell^p(\bZ)$ is the operator $W$ given by $W e_k=e_{k+1}$ for all $k\in \bZ$.
It is clear that $W$ is an invertible isometry such that $W\eN_k=\eN_{k+1}$ and $W^{-1}\eN_k=\eN_{k-1}$,
for all $k\in \bZ$. Hence, the group $\pO(\fN)=\{ W^k;\;k\in \bZ\}$ is contained in $\Col(\fN)$.

\begin{proposition} \label{prop12}
$\Col(\fN)=\Grp\bigl(\Alg(\fN)\bigr)\rtimes \pO(\fN)$.
\end{proposition}

\begin{proof}
It is obvious that $\pO(\fN)\cap \Grp\bigl(\Alg(\fN)\bigr)=\{ I\}$. Let $S\in \Col(\fN)$ be arbitrary. Then there
exists an integer $s$ such that $S\eN_0=\eN_{s}$. For an integer $k\geq 0$, we have
$\eN_{-k}\subsetneq \cdots\subsetneq \eN_0\subsetneq \cdots\subsetneq \eN_k$
and therefore $S\eN_{-k}\subsetneq\cdots\subsetneq \eN_s=S\eN_0\subsetneq \cdots\subsetneq S\eN_k$.
Since $S$ is invertible, the codimension of $\eN_s$ in $S\eN_k$ is $k$, so that $\eN_{s+k}=S\eN_k$. Similarly one can see that 
$\eN_{s-k}=S\eN_{-k}$. To finish the proof, observe that we can write $S=(SW^{-k})W^k\in \Grp\bigl(\Alg(\fN)\bigr)\pO(\fN)$.
\end{proof}

We now consider a family of nests in the Banach space $C_0[-1,1]$ of continuous complex-valued functions
on the interval $[-1,1]$ that vanish at the endpoints. Let $\varBeta([-1,1])$ be the family of all
monotone bijections $\psi\colon[-1,1]\to[-1,1]$. Of course, $\varBeta([-1,1])$ is a group under the composition.
Since every monotone function $\psi\colon [a,b]\to[c,d]$ is continuous on $[a,b]$ except at possibly countably 
many jump points (see \cite[Theorem 1, \S 6.1]{RF}), it follows that every monotone bijection 
$\psi\colon [a,b]\to[c,d]$ is continuous. Therefore, the group $B([-1,1])$ consists of continuous functions with 
the identity function $\iota\colon x\mapsto x$ as its unit element. Let $\varGamma([-1,1])$ be the subgroup
consisting of all increasing functions in $\varBeta([-1,1])$. Notice that $\varBeta([-1,1])$ is a semidirect product of
$\varGamma([-1,1])$ and the subgroup $\{ \iota, -\iota\}$ of $\varBeta([-1,1])$. In particular, every decreasing
function in $\varBeta([-1,1])$ is of the form $x\mapsto \psi(-x)$ for some function $\psi\in\varGamma([-1,1])$.

Recall that the support of $f\in C_0[-1,1]$ is $\supp(f)=\overline{\{x\in[-1,1];\; f(x)\ne 0\}}$. Let
$\theta_1\colon[0,1]\to[-1,0]$ be a decreasing bijective function, and let $\theta_2\colon[0,1]\to[0,1]$ be
an increasing bijective function. By the observation above, $\theta_1$ and $\theta_2$ are continuous.
For each $t\in[0,1]$, let $\eN_t=\bigl\{f\in C_0[-1,1];\; \supp(f)\subseteq[\theta_1(t),\theta_2(t)]\bigr\}$.
It is clear that every $\eN_t$ is a subspace of $C_0[-1,1]$. Notice that $\eN_0=\{ 0\}$ and $\eN_1=C_0[-1,1]$;
moreover, $\eN_{t_1}\subsetneq \eN_{t_2}$ whenever $0\leq t_1 < t_2\leq 1$. It is not hard to see that
$\fN=\{ \eN_{t};\; t\in[0,1]\}$ is a maximal nest in $\fC\bigl(C_0[-1,1]\bigr)$.

For every $\varphi\in\varGamma([-1,1])$, the mappings $V_\varphi\colon f\mapsto f\circ\varphi$ and
$V_{\varphi^{-1}}\colon f\mapsto f\circ\varphi^{-1}$ are well defined on $C_0[-1,1]$. It is easily seen that they
are isometric isomorphisms such that $V_{\varphi}^{-1}=V_{\varphi^{-1}}$.

\begin{theorem} \label{theo09}
The following statements hold. 
\begin{itemize}
    \item [(i)] $\varGamma_{\theta_1,\theta_2}([-1,1]):=\{ \varphi\in \varGamma([-1,1]);\;
\theta_{1}^{-1}\circ\varphi\circ\theta_1=\theta_{2}^{-1}\circ\varphi\circ\theta_2\}$
is a subgroup of $\varGamma([-1,1])$.
    \item [(ii)] If $\varphi\in\varGamma([-1,1])$, then $V_\varphi$ is a collineation of $\fN$ if and only if
$\varphi\in\varGamma_{\theta_1,\theta_2}([-1,1])$. Hence, $\pO(\fN)=\{ V_\varphi;\; \varphi\in
\varGamma_{\theta_1,\theta_2}([-1,1])\}$ is a subgroup of $\Col(\fN)$.
    \item [(iii)] $\Col(\fN)=\Grp\bigl(\Alg(\fN)\bigr)\rtimes \pO(\fN)$.
\end{itemize} 
\end{theorem}

\begin{proof}
(i) If $\varphi\in\varGamma_{\theta_1,\theta_2}([-1,1])$, then
$\theta_{1}^{-1}\bigl(\varphi(0)\bigr)=\theta_{2}^{-1}\bigl(\varphi(0)\bigr)$. This implies that $\varphi(0)$ lies 
in the intersection of the domains of $\theta_{1}^{-1}$ and $\theta_{2}^{-1}$, and therefore $\varphi(0)=0$. 
It follows that $\varphi$ is an increasing bijection of both $[-1,0]$ and of $[0,1]$. Consequently, the map
$\theta_{1}^{-1}\circ\varphi\circ\theta_1=\theta_{2}^{-1}\circ\varphi\circ\theta_2$ is a bijection of $[0,1]$.
Its inverse is given by
$\theta_{1}^{-1}\circ\varphi^{-1}\circ\theta_1=\theta_{2}^{-1}\circ\varphi^{-1}\circ\theta_2$, which implies that
$\varphi^{-1}\in\varGamma_{\theta_1,\theta_2}([-1,1])$. It is straightforward to see that the composition of two
functions from $\varGamma_{\theta_1,\theta_2}([-1,1])$ is in $\varGamma_{\theta_1,\theta_2}([-1,1])$, as well.

(ii) By \Cref{prop15}, $V_\varphi\in\Col(\fN)$ if and only if $V_\varphi\eN_t\in \fN$ for every $t\in[0,1]$.
Assume that $\varphi\in\varGamma_{\theta_1,\theta_2}([-1,1])$ and let $t\in[0,1]$ be arbitrary.
If $f\in \eN_t$, then $\supp(f)\subseteq[\theta_1(t),\theta_2(t)]$ and therefore
$\supp(f\circ\varphi)\subseteq\bigl[ \varphi^{-1}(\theta_1(t)),\varphi^{-1}(\theta_2(t))\bigr]$.
Denote $s=\bigl(\theta_{1}^{-1}\circ\varphi^{-1}\circ\theta_1\bigr)(t)$. Since
$\varphi^{-1}\in \varGamma_{\theta_1,\theta_2}([-1,1])$ we also have
$s=\bigl(\theta_{2}^{-1}\circ\varphi^{-1}\circ\theta_2\bigr)(t)$. Hence,
$\varphi^{-1}(\theta_1(t))=\theta_1(s)$ and $\varphi^{-1}(\theta_2(t))=\theta_2(s)$. This shows that
$\bigl[ \varphi^{-1}(\theta_1(t)),\varphi^{-1}(\theta_2(t))\bigr]=\bigl[ \theta_1(s),\theta_2(s)\bigr]$.
We have proved that $V_\varphi f\in \eN_s$ whenever $f\in \eN_t$. On the other hand, a similar reasoning
shows that $V_{\varphi^{-1}}g\in \eN_t$ if $g\in\eN_s$. This proves that $V_\varphi \eN_t=\eN_s$, that is,
$V_\varphi\in\Col(\fN)$.

Suppose that $\varphi\in\varGamma([-1,1])$ and $V_\varphi\in\Col(\fN)$. Hence, for every $t\in[0,1]$
there exists a unique $s\in[0,1]$ such that $V_\varphi\eN_t=\eN_s$. Since
$\supp(f\circ\varphi)=\varphi^{-1}\bigl( \supp(f)\bigr)$ for every $f\in C_0[-1,1]$, we see that
$[\theta_1(s),\theta_2(s)]=\varphi^{-1}\bigl([\theta_1(t),\theta_2(t)]\bigr)$. It follows that
$\theta_1(s)=\varphi^{-1}\bigl(\theta_1(t)\bigr)$ and $\theta_2(s)=\varphi^{-1}\bigl(\theta_2(t)\bigr)$
and consequently $\varphi\in\varGamma_{\theta_1,\theta_2}([-1,1])$. Moreover,
\begin{equation} \label{eq01}
V_\varphi \eN_t=\eN_{(\theta_{1}^{-1}\circ\varphi^{-1}\circ\theta_1)(t)}=
\eN_{(\theta_{2}^{-1}\circ\varphi^{-1}\circ\theta_2)(t)},\qquad \text{for every}\; t\in[0,1].
\end{equation}

It is not hard to see that $\pO(\fN)$ is a subgroup of $\Col(\fN)$ (notice that $\varphi\mapsto V_\varphi$
is an injective anti-homomorphism from $\varphi\in\varGamma([-1,1])$ into $\Col(\fN)$).

(iii) If $V_\varphi\in \Grp\bigl(\Alg(\fN)\bigr)$, then
$\varphi^{-1}\bigl([\theta_1(t),\theta_2(t)]\bigr)=[\theta_1(t),\theta_2(t)]$ for every $t\in[0,1]$
(see the proof of (ii)). It follows that $\varphi\bigl(\theta_1(t)\bigr)=\theta_1(t)$ and
$\varphi\bigl(\theta_2(t)\bigr)=\theta_2(t)$ for every $t\in[0,1]$. Since $\theta_1\colon[0,1]\to[-1,0]$
and $\theta_2\colon[0,1]\to[0,1]$ are bijections equality $\varphi\bigl(\theta_1(t)\bigr)=\theta_1(t)$ gives
$\varphi(x)=x$ for every $x\in[-1,0]$, and equality $\varphi\bigl(\theta_2(t)\bigr)=\theta_2(t)$ gives
$\varphi(x)=x$ for every $x\in[0,1]$. Hence, $\varphi$ is the identity function on $[-1,1]$ and therefore
$V_\varphi$ is the identity operator. This proves that $\Grp\bigl(\Alg(\fN)\bigr)\cap \pO(\fN)=\{ I\}$.

Let $S\in \Col(\fN)$ be an arbitrary collineation. For every $t\in[0,1]$, there exists a unique $s\in[0,1]$
such that $S\eN_t=\eN_s$. Hence, there is a correspondence
$[\theta_1(t),\theta_2(t)]\leftrightarrow[\theta_1(s),\theta_2(s)]$, more precisely, there are functions
$\mu\colon[0,1]\to[0,1]$ and $\nu\colon[-1,0]\to[-1,0]$ such that
$[\theta_1(s),\theta_2(s)]=[\nu(\theta_1(t)),\mu(\theta_2(t))]$. It follows that
$\theta_1(s)=\nu(\theta_1(t))$ and $\theta_2(s)=\mu(\theta_2(t))$ and therefore
$\bigl(\theta_{1}^{-1}\circ\nu\circ\theta_1\bigr)(t)=\bigl(\theta_{2}^{-1}\circ\mu\circ\theta_2\bigr)(t)$,
for every $t\in[0,1]$.
Since $\eN_t\mapsto S\eN_t$ is an increasing bijection, we see that $\mu$ and $\nu$ are increasing bijections, too.
Hence, they are continuous functions such that $\mu(0)=0=\nu(0)$. Let $\varphi\colon[-1,1]\to[-1,1]$
be given by
$$\varphi(x)=\left\{\begin{array}{cr}
\nu(x),&\; x\in[-1,0),\\
\mu(x),&\; x\in[0,1].
\end{array}\right. $$
Then $\varphi$ is a continuous increasing bijection such that
$\theta_{1}^{-1}\circ\varphi\circ\theta_1=\theta_{2}^{-1}\circ\varphi\circ\theta_2$. Hence,
$\varphi\in \varGamma_{\theta_1,\theta_2}([-1,1])$ and therefore $V_\varphi\in\Col(\fN)$, by (ii).
Notice that $S\eN_t=\eN_{(\theta_{1}^{-1}\circ\varphi^{-1}\circ\theta_1)(t)}=
\eN_{(\theta_{2}^{-1}\circ\varphi^{-1}\circ\theta_2)(t)}$ for every $t\in[0,1]$. It follows that (see \eqref{eq01})
$$SV_{\varphi}^{-1}\eN_t=S\bigl(\eN_{(\theta_{2}^{-1}\circ\varphi\circ\theta_2)(t)}\bigr)=
\eN_{(\theta_{2}^{-1}\circ\varphi^{-1}\circ\theta_2\circ\theta_{2}^{-1}\circ\varphi\circ\theta_2)(t)}=\eN_t,
\qquad (t\in[0,1]).$$
We see that $S=(SV_{\varphi}^{-1}) V_{\varphi}$, where $SV_{\varphi}^{-1}\in\Grp\bigl(\Alg(\fN)\bigr)$
and $V_{\varphi}\in \pO(\fN)$.
\end{proof}

Let $\fN\subseteq \fC\bigl(C_0[-1,1]\bigr)$ be as before and let $\fN_\bQ\subseteq \fN$ be the family of
all subspaces $\eN_t$ with $t\in[0,1]\cap\bQ$. It is clear that $\fN_\bQ$ is a totally ordered lattice
of subspaces. It is not complete, its completion is $\fN$, that is, $\widehat{\fN_\bQ}=\fN$ in the notation
we have introduced in the paragraph preceding \Cref{theo07}. An operator $V_\varphi\in \pO(\fN)$
is a collineation of $\fN_\bQ$ if and only if function $\theta_{2}^{-1}\circ\varphi^{-1}\circ\theta_2\colon[0,1]\to
[0,1]$ maps rational points to rational points. Of course, not every function
$\varphi\in \varGamma_{\theta_1,\theta_2}([-1,1])$ has this property. For instance, let $\theta_1(t)=-t$
and $\theta_2(t)=t$ for $t\in[0,1]$. In this case, an increasing bijection $\varphi\colon[-1,1]\to[-1,1]$
belongs to $\varGamma_{\theta_1,\theta_2}([-1,1])$ if and only if $\varphi(-x)=-\varphi(x)$ for every $x\in[-1,1]$.
In particular, the function $\varphi(x)=x^3$ is in $\varGamma_{\theta_1,\theta_2}([-1,1])$. Since
$\bigl(\theta_{2}^{-1}\circ\varphi^{-1}\circ\theta_2\bigr)(t)=\varphi^{-1}(t)=\sqrt[3]{t}$ $(t\in[0,1])$
does not map rational points to rational points we see that the corresponding isometry $V_\varphi$ is not 
a collineation of $\fN_\bQ$. This shows that the group of collineations of a lattice of subspaces $\fL$ can be a 
proper subgroup of the group of collineations of the completion $\widehat{\fL}$ (cf. \Cref{theo07}).

\subsection{Volterra nests in $L^p [0,1]$} \label{Sec053}
Let $m$ be the Lebesgue measure on $[0,1]$ and let $L^0[0,1]$ be the vector space of all equivalence classes of 
measurable functions on $[0,1]$. As is customary, we will abuse notation slightly by writing $f$ in place of 
the equivalence class $[f]$. For $1\leq p<\infty$, let $L^p [0,1]\subseteq L^0[0,1]$ be the Banach space of 
Lebesgue $p$-integrable functions on $[0,1]$ (see \cite[Ch. 3]{Rud} for details).
For a measurable set $E\subseteq [0,1]$, let $\chi_E$ be the characteristic function of $E$. It is obvious that
$\chi_E\in L^p[0,1]$ for every $1\leq p< \infty$. For each $t\in[0,1]$, let
$\eN_{t}=\{ f\in L^p[0,1];\; f\chi_{[t,1]}=0\; a.e.\}$. Clearly, $\eN_{t}$ is a subspace of $L^p [0,1]$.
Notice that $\eN_{0}=\{0\}$, $\eN_{1}=L^p [0,1]$, and $\eN_{t_1}\subsetneq \eN_{t_2}$
if $0\leq t_1<t_2\leq 1$. The family $\fN=\{ \eN_{t};\, t\in[0,1]\}$ is a (continuous) nest in $\fC(L^p [0,1])$, 
called the Volterra nest.

Let $\varGamma$ be the group of all increasing bijective functions on $[0,1]$. As already mentioned, 
$\varGamma$ consists of continuous functions. By Lebesgue's Theorem \cite[\S 6.2]{RF}, every function
$\varphi\in \varGamma$ is differentiable almost everywhere on the interval $(0,1)$. Since $\varphi$ is increasing, 
we have $\varphi'\geq 0$ on the set where $\varphi$ is differentiable. Moreover, since $\varphi^{-1}$ is continuous,  
$\varphi$ maps measurable sets to measurable sets. It follows that
$m_\varphi(E)=m\bigl(\varphi(E)\bigr)$ defines a measure $m_\varphi$ on $[0,1]$.

Recall that $\psi\colon [0,1]\to \bR$ is an absolutely continuous function if
for each $\varepsilon>0$ there exists $\delta(\varepsilon)>0$ such that for each $n\in \bN$ and numbers
$0\leq x_1\leq y_1\leq \cdots \leq x_n\leq y_n\leq 1$ such that $\sum_{i=1}^{n}(y_i-x_i)<\delta(\varepsilon)$
we have $\sum_{i=1}^{n}|\psi(y_i) - \psi(x_i)| < \varepsilon$. 
The proof of the following result can be found in \cite[Theorem 7.18]{Rud}.

\begin{theorem}\label{theo03}
If $\varphi\in\varGamma$, then the following assertions are equivalent: 
\begin{itemize}
    \item [(i)] $\varphi$ is absolutely continuous;
    \item [(ii)] the derivative $\varphi'$ is in $L^1[0,1]$ and $\varphi(x)=\int_{0}^{x}\varphi'(t)\, dt$ 
    for any $x\in[0,1]$;
    \item [(iii)] the measure $m_\varphi$ is absolutely continuous with respect to the measure $m$.\qed
\end{itemize}
\end{theorem}

The following example shows that there exist absolutely continuous functions in $\varGamma$
whose inverses are not absolutely continuous.

\begin{example} \label{ex03}
This example is based on a comment in the last paragraph of \S 3.3 in \cite{KS}.
Let $c\colon[0,1]\to[0,1]$ be the usual Cantor-Lebesgue function (see \cite[\S 2.7]{RF}) and let
$\psi(x)=\frac{1}{2}(x+c(x))$ for every $x\in[0,1]$. It is easily seen that $\psi\in \varGamma$. This function is
not absolutely continuous since it maps the Cantor set, which has measure $0$, onto a measurable set of
positive measure  (see \cite[\S 2.7 Proposition 21]{RF}).
On the other hand, the inverse $\varphi=\psi^{-1}$ is absolutely continuous, and even more it is a Lipschitz function.
Indeed, if $x_1, x_2\in[0,1]$ and $x_1< x_2$, then $\varphi(x_1)<\varphi(x_2)$ because $\varphi$ is an
increasing bijection. Since $x=\psi(\varphi(x))=\frac{1}{2}\bigl(\varphi(x)+c(\varphi(x))\bigr)$ for each $x\in [0,1]$, 
we have
$$\left|\frac{x_2-x_1}{\varphi(x_2)-\varphi(x_1)}\right|
=\frac{\varphi(x_2)+c(\varphi(x_2))-\varphi(x_1)-c(\varphi(x_1))}{2(\varphi(x_2)-\varphi(x_1))}
=\frac{1}{2}+\frac{c(\varphi(x_2))-c(\varphi(x_1))}{2(\varphi(x_2)-\varphi(x_1))}\geq \frac{1}{2}
$$
as $c$ is non-decreasing. Hence, $\varphi$ is a Lipschitz function and therefore absolutely continuous, 
by \cite[\S 6.4 Proposition 7]{RF}, but the inverse of $\varphi$ is not absolutely continuous. \qed
\end{example}

The next lemma is known. For the sake of completeness, we include a short proof.

\begin{lemma} \label{lem04}
The following assertions are equivalent for a function $\varphi\in\varGamma$:

\begin{itemize}
\item[(i)] $\varphi$ and $\varphi^{-1}$ are absolutely continuous;
\item[(ii)] measures $m$ and $m_\varphi$ are mutually absolutely continuous;
\item[(iii)] $\varphi$ is absolutely continuous and $\varphi'>0$ a.e.\
\end{itemize}
\end{lemma}

\begin{proof}
(i)$\iff$(ii). If $\varphi$ and $\varphi^{-1}$ are absolutely continuous, by \Cref{theo03} the measures $m_\varphi$
and $m_{\varphi^{-1}}$ are absolutely continuous with respect to the Lebesgue measure $m$. To see that
$m$ is absolutely continuous with respect to $m_\varphi$, choose a measurable set $E\subseteq [0,1]$ such that
$m_\varphi(E)=0$. Hence, $m\bigl(\varphi(E)\bigr)=0$ and therefore $m(E)=
m_{\varphi^{-1}}\bigl(\varphi(E)\bigr)=0$ since $m_{\varphi^{-1}}$ is absolutely continuous with respect to $m$.

To prove the opposite implication, suppose that $m$ and $m_\varphi$ are mutually absolutely continuous.
Let $E\subseteq[0,1]$ be a measurable set such that $m(E)=0$. Since
$m(E)=m_{\varphi}\bigl(\varphi^{-1}(E)\bigr)$ and $m$ is absolutely continuous with respect to $m_\varphi$
we have $m\bigl(\varphi^{-1}(E)\bigr)=0$. This proves that $m_{\varphi^{-1}}\ll m$. The proof of
$m_{\varphi}\ll m$ is similar.

(i)$\iff$(iii). Assume first that $\varphi$ and $\varphi^{-1}$ are absolutely continuous. Let $D_\varphi$ and
$D_{\varphi^{-1}}$ be the subsets of $[0,1]$ where $\varphi$ and $\varphi^{-1}$, respectively, are differentiable.  
By Lebegue's theorem, the complements of $D_\varphi$ and $D_{\varphi^{-1}}$ have Lebesgue measure zero. 
Since $\varphi^{-1}$ is absolutely continuous and bijective, the set
$\varphi^{-1}([0,1]\setminus D_{\varphi^{-1}})=[0,1]\setminus\varphi^{-1}(D_{\varphi^{-1}})$ has 
Lebesgue measure zero as well.
It follows that the complement of $\varphi^{-1}(D_{\varphi^{-1}})\cap D_{\varphi}$ has Lebesgue measure zero.
For $x\in \varphi^{-1}(D_{\varphi^{-1}})\cap D_{\varphi}$, the derivatives $(\varphi^{-1})'(\varphi(x))$ and
$\varphi'(x)$ exist and are non-negative. Since $\varphi^{-1}(\varphi(x))=x$ we have
$(\varphi^{-1})'(\varphi(x)) \varphi'(x)=1$, and therefore, $\varphi'(x)>0$. This proves that
$\varphi'>0$ a.e.

Suppose now that $\varphi$ is absolutely continuous and $\varphi'>0$ a.e.\ Denote
$P_\varphi=\{ x\in D_\varphi; \varphi'(x)>0\}$. By assumption, $m\bigl([0,1]\setminus P_\varphi\bigr)=0$.
Hence, $m\bigl( \varphi^{-1}(E)\cap([0,1]\setminus P_\varphi)\bigr)=0$ for every measurable set
$E\subseteq [0,1]$.
By \cite[Lemma 5.8.13]{Bog}, it follows that  if $m(E)=0$, then 
$m\bigl(\varphi^{-1}(E)\bigr)=m\bigl(\varphi^{-1}(E)\cap P_\varphi\bigr)=0.$
Hence, $\varphi^{-1}$ is absolutely continuous.
\end{proof}

Notice that \Cref{lem04} and \Cref{ex03} show that there exists an absolutely continuous function $\varphi$
such that set $\{ x\in D_\varphi; \varphi'(x)=0\}$ has positive measure. Hence, the claim {\em ``if $\varphi$
is absolutely continuous then $\varphi'>0$ a.e.''} on page 26 of \cite{AAW} is not correct.

We now return to the study of collineations of the Volterra nest $\fN$. The first observation is simple.
If $\varphi\in\varGamma$, then the composition $f\mapsto f\circ\varphi$ induces a well-defined linear transformation
on $L^0[0,1]$ if and only if $\varphi^{-1}$ is absolutely continuous. To see this, assume first that $\varphi^{-1}$ 
is absolutely continuous. If measurable functions $f$ and $g$ represent the same equivalence class in $L^0[0,1]$, 
then $\{ x\in[0,1];\; \bigl(f\circ\varphi\bigr)(x)\ne \bigl(g\circ\varphi\bigr)(x)\}=
\varphi^{-1}\bigl(\{x\in[0,1];\; f(x)\ne g(x)\}\bigr)$
is a set of measure $0$ since $\{x\in[0,1];\; f(x)\ne g(x)\}$ has measure $0$ and $\varphi^{-1}$ maps sets of
measure $0$ to sets of measure $0$. Thus, $f\mapsto f\circ\varphi$ induces a well-defined transformation on 
$L^0[0,1]$. On the other hand, if the composition  $f\mapsto f\circ\varphi$ induces a well-defined transformation 
on $L^0[0,1]$, then $\chi_E\circ\varphi=0$ a.e. for every measurable set $E\subseteq[0,1]$ with measure $0$. 
Since $\chi_E\circ\varphi=\chi_{\varphi^{-1}(E)}$ we see that $\varphi^{-1}(E)$ has measure $0$ whenever $E$ 
has measure $0$. We conclude that $\varphi^{-1}$ is absolutely continuous.

Denote by $\varDelta$ the family of all functions $\varphi\in\varGamma$ satisfying the equivalent conditions 
of \Cref{lem04}. It is easily seen that $\varDelta$ is a subgroup of $\varGamma$.

\begin{proposition} \label{prop03}
Let $1\leq p< \infty$. If $\varphi\in \varDelta$, then 
$V_{\varphi}\colon f\mapsto (\varphi')^{1/p} \bigl(f\circ\varphi\bigr)$ is an isometric operator on $L^p[0,1]$ with
the inverse $V_{\varphi^{-1}}$. The mapping $\omega\colon \varphi\mapsto V_{\varphi}$ is an injective 
anti-homomorphism from $\varDelta$ into the group $\Col(\fN)$.
\end{proposition}

\begin{proof}
By the above discussion, if $\varphi\in\varDelta$, then $f\mapsto f\circ\varphi$ and $f\mapsto f\circ\varphi^{-1}$ are
well-defined linear transformations on $L^0[0,1]$. Furthermore, \Cref{lem04} yields that 
$\varphi'>0$ a.e.\ and $(\varphi^{-1})'>0$ a.e., so that $V_{\varphi}$ and $V_{\varphi^{-1}}$ are well-defined
linear transformations on $L^0[0,1]$.

We claim that $V_\varphi$ and $V_{\varphi^{-1}}$ are isometric isomorphisms of $L^p[0,1]$ 
and that $V_{\varphi^{-1}}=V_\varphi^{-1}$. Since $\varphi\colon[0,1]\to[0,1]$ is an increasing absolutely
continuous bijection, by the special case of the change-of-variables theorem \cite[Theorem 7.26]{Rud}, we have 
$$ \|V_{\varphi} f\|_{p}^{p}= \int_{0}^{1} \varphi'(x) |\bigl(f\circ\varphi\bigr)(x)|^p\, dx=
\int_{0}^{1} |f(x)|^p\, dx=\|f\|_{p}^{p}\qquad (f\in L^p[0,1]),$$
proving that $V_\varphi\colon L^p[0,1]\to L^p[0,1]$ is an isometry. A similar argument yields that 
$V_{\varphi^{-1}}\colon L^p[0,1]\to L^p[0,1]$ is an isometry. 
A direct calculation shows $V_{\varphi}^{-1}=V_{\varphi^{-1}}$.

Suppose $f\in \eN_t$. Then $V_\varphi f\in \eN_{\varphi^{-1}(t)}$, since the condition $f(x)=0$ for almost every 
$x\in [t,1]$ implies that $f(\varphi(x))=0$ for almost every $x\in [\varphi^{-1}(t),1]$. On the other hand, since 
$\varphi'>0$ a.e.\ we see that $V_{\varphi}\chi_{[0,t]}$ is not contained in $\eN_{s}$ for any $s<\varphi^{-1}(t)$. 
Thus, $V_{\varphi}\eN_{t}=\eN_{\varphi^{-1}(t)}$. Similarly, $V_{\varphi^{-1}}\eN_{t}=\eN_{\varphi(t)}$. 
This shows that $V_{\varphi}$ is a collineation of $\fN$. It is straightforward to check that 
$V_{\varphi\circ\psi}=V_{\psi}V_{\varphi}$ which means $\omega$ is an anti-homomorphism. It is injective 
since $V_{\varphi}=V_{\psi}$ gives $V_{\varphi\circ\psi^{-1}}=I$ and therefore $(\varphi\circ\psi^{-1})(x)=x$ 
for every $x\in[0,1]$, that is, $\varphi=\psi$.
\end{proof}

Denote by $\pO(\fN)$ the subgroup of isometries in $\Col(\fN)$ induced by functions in $\varDelta$, that is,
$\pO(\fN)=\omega(\varDelta)$. 

\begin{theorem} \label{theo08}
For every $1\leq p< \infty$, the semidirect product $\Grp\bigl(\Alg(\fN)\bigr)\rtimes \pO(\fN)$ is a subgroup of 
$\Col(\fN)$.

\begin{itemize}
\item[(i)] If $p=1$, then $\Col(\fN)=\Grp\bigl(\Alg(\fN)\bigr)\rtimes \pO(\fN)$.
\item [(ii)] Let $p=2$. For every $\varphi\in \varGamma$ there exists $S\in\Col(\fN)$ such that
$S\eN_{t}=\eN_{\varphi(t)}$ $(t\in[0,1])$. On the other hand, there exists a unitary operator $U\in\Col(\fN)$ 
such that $U\eN_{t}=\eN_{\varphi(t)}$ $(t\in[0,1])$ if and only if $\varphi\in \varDelta$. Hence, 
if $\pO\subseteq \Col(\fN)$ is any subgroup of unitary operators, then $\Grp\bigl(\Alg(\fN)\bigr) \pO$
is a proper subgroup of $\Col(\fN)$; in particular this holds for the semidirect product
$\Grp\bigl(\Alg(\fN)\bigr)\rtimes \pO(\fN)$.
\end{itemize}
\end{theorem}

\begin{proof} It is clear that $\Grp\bigl(\Alg(\fN)\bigr)\cap\pO(\fN)=\{I\}$ and therefore the semidirect
product $\Grp\bigl(\Alg(\fN)\bigr)\rtimes \pO(\fN)$ is a subgroup of $\Col(\fN)$.

(i) Let $S\in\Col(\fN)$ be arbitrary. It is not hard to see that there is a unique function $\varphi\in \varGamma$
such that $S\eN_t=\eN_{\varphi(t)}$ and therefore $S^{-1}\eN_t=\eN_{\varphi^{-1}(t)}$ for every $t\in[0,1]$.
Define $\eM_t(\varphi):=\eN_{\varphi(t)}$ and notice that $\fM(\varphi)=\{ \eM_t(\varphi); t\in[0,1]\}$ 
is a continuous nest consisting of the same subspaces as $\fN$ with a possibly different indexing. 
Hence, $S$ is an invertible operator that implements the order isomorphism 
$\theta_\varphi\colon \eN_t\mapsto \eM_t(\varphi)$. By \cite[Theorem 4.4]{ALWW}, this can happen if and only if 
$\varphi$ and $\varphi^{-1}$ are absolutely continuous functions.
Since $SV_\varphi\in \Grp\bigl(\Alg(\fN)\bigr)$ and $S=(SV_\varphi)V_{\varphi{-1}}$ we conclude that
$S\in \Grp\bigl(\Alg(\fN)\bigr)\rtimes \pO(\fN)$.

(ii) For $\varphi\in\varGamma$, let $\eM_t(\varphi)$ $(t\in [0,1])$ and $\fM(\varphi)$
have the same meaning as in the proof of (i). By \cite[Theorem 3.2]{Dav1}, there exists an invertible operator
$S\in \pB(L^2[0,1])$ such that the order isomorphism $\theta_\varphi\colon \eN_t\mapsto \eM_t(\varphi)$ 
is implemented by $S$, that is, $S\eN_t=\eM_t(\varphi)=\eN_{\varphi(t)}$ for every $t\in[0,1]$. 
Since this gives $S^{-1}\eN_t=\eN_{\varphi^{-1}(t)}$ we conclude that $S$ is a collineation of $\fN$.
It follows from \cite[\S 7]{Dav} (see Example 7.18 in that chapter) that $U\eN_{t}=\eN_{\varphi(t)}$ $(t\in[0,1])$
for a unitary operator $U$ if and only if $\varphi\in \varDelta$. The first proof of this fact is probably given by
Kadison and Singer (see Theorem 3.3.1 in \cite{KS} and the paragraph after its proof).

Let $\pO\subseteq \Col(\fN)$ be a subgroup of unitary operators. Suppose that $\varphi\in \varGamma$ is 
not in $\varDelta$. We claim that the corresponding collineation $S$ is not in the subgroup 
$\Grp\bigl(\Alg(\fN)\bigr) \pO$. Notice that $\Grp\bigl(\Alg(\fN)\bigr) \pO=\pO\Grp\bigl(\Alg(\fN)\bigr)$ since
$\Grp\bigl(\Alg(\fN)\bigr)$ is a normal subgroup of $\Col(\fN)$. To prove our claim, assume that $S$ were in 
$\pO\Grp\bigl(\Alg(\fN)\bigr)$. Then there would exist $U\in\pO$ and $A\in\Grp\bigl(\Alg(\fN)\bigr)$ such that
$S=UA$. Hence, we would have $\eN_{\varphi(t)}=S\eN_t=UA\eN_t=U\eN_t$ for every $t\in[0,1]$ which 
is impossible since $\varphi\not\in\varDelta$. Thus, $\Grp\bigl(\Alg(\fN)\bigr) \pO$, in particular 
$\Grp\bigl(\Alg(\fN)\bigr)\rtimes \pO(\fN)$, is a proper subgroup of $\Col(\fN)$.
\end{proof}

\Cref{theo08} partially answers the question of whether $\Grp\bigl(\Alg(\fN)\bigr)$ is complemented in $\Col(\fN)$
in the case when $\fN$ is the Volterra nest in $\fC\bigl(L^p[0,1]\bigr)$.
The problem is completely solved for $p=1$, but for $1<p<\infty$ it remains open. What we know more
is that in the case $p=2$ the group $\Grp\bigl(\Alg(\fN)\bigr)$ cannot be complemented by a group of
isometries (i.e., unitary operators).

\subsection*{Acknowledgements}
The first author was supported by the Slovenian Research and Innovation Agency program P2-0268.
The second author was supported by the Slovenian Research and Innovation Agency program P1-0222 
and grant J1-50002.
%
%
	
\end{document}